
\font\ninerm=cmr9  \font\eightrm=cmr8  \font\sixrm=cmr6
\font\ninei=cmmi9  \font\eighti=cmmi8  \font\sixi=cmmi6
\font\ninesy=cmsy9 \font\eightsy=cmsy8 \font\sixsy=cmsy6
\font\ninebf=cmbx9 \font\eightbf=cmbx8 \font\sixbf=cmbx6
\font\nineit=cmti9 \font\eightit=cmti8 
\font\ninett=cmtt9 \font\eighttt=cmtt8 
\font\ninesl=cmsl9 \font\eightsl=cmsl8

\font\twelverm=cmr12 at 15pt
\font\twelvei=cmmi12 at 15pt
\font\twelvesy=cmsy10 at 15pt
\font\twelvebf=cmbx12 at 15pt
\font\twelveit=cmti12 at 15pt
\font\twelvett=cmtt12 at 15pt
\font\twelvesl=cmsl12 at 15pt
\font\twelvegoth=eufm10 at 15pt

\font\tengoth=eufm10  \font\ninegoth=eufm9
\font\eightgoth=eufm8 \font\sevengoth=eufm7 
\font\sixgoth=eufm6   \font\fivegoth=eufm5
\newfam\gothfam \def\goth{\fam\gothfam\tengoth} 
\textfont\gothfam=\tengoth
\scriptfont\gothfam=\sevengoth 
\scriptscriptfont\gothfam=\fivegoth

\catcode`@=11
\newskip\ttglue

\def\tenpoint{\def\rm{\fam0\tenrm}
  \textfont0=\tenrm \scriptfont0=\sevenrm
  \scriptscriptfont0\fiverm
  \textfont1=\teni \scriptfont1=\seveni
  \scriptscriptfont1\fivei 
  \textfont2=\tensy \scriptfont2=\sevensy
  \scriptscriptfont2\fivesy 
  \textfont3=\tenex \scriptfont3=\tenex
  \scriptscriptfont3\tenex 
  \textfont\itfam=\tenit\def\it{\fam\itfam\tenit}%
  \textfont\slfam=\tensl\def\sl{\fam\slfam\tensl}%
  \textfont\ttfam=\tentt\def\tt{\fam\ttfam\tentt}%
  \textfont\gothfam=\tengoth\scriptfont\gothfam=\sevengoth 
  \scriptscriptfont\gothfam=\fivegoth
  \def\goth{\fam\gothfam\tengoth}
  \textfont\bffam=\tenbf\scriptfont\bffam=\sevenbf
  \scriptscriptfont\bffam=\fivebf
  \def\bf{\fam\bffam\tenbf}%
  \tt\ttglue=.5em plus.25em minus.15em
  \normalbaselineskip=12pt \setbox\strutbox\hbox{\vrule
  height8.5pt depth3.5pt width0pt}%
  \let\big=\tenbig\normalbaselines\rm}

\def\ninepoint{\def\rm{\fam0\ninerm}
  \textfont0=\ninerm \scriptfont0=\sixrm
  \scriptscriptfont0\fiverm
  \textfont1=\ninei \scriptfont1=\sixi
  \scriptscriptfont1\fivei 
  \textfont2=\ninesy \scriptfont2=\sixsy
  \scriptscriptfont2\fivesy 
  \textfont3=\tenex \scriptfont3=\tenex
  \scriptscriptfont3\tenex 
  \textfont\itfam=\nineit\def\it{\fam\itfam\nineit}%
  \textfont\slfam=\ninesl\def\sl{\fam\slfam\ninesl}%
  \textfont\ttfam=\ninett\def\tt{\fam\ttfam\ninett}%
  \textfont\gothfam=\ninegoth\scriptfont\gothfam=\sixgoth 
  \scriptscriptfont\gothfam=\fivegoth
  \def\goth{\fam\gothfam\tengoth}
  \textfont\bffam=\ninebf\scriptfont\bffam=\sixbf
  \scriptscriptfont\bffam=\fivebf
  \def\bf{\fam\bffam\ninebf}%
  \tt\ttglue=.5em plus.25em minus.15em
  \normalbaselineskip=11pt \setbox\strutbox\hbox{\vrule
  height8pt depth3pt width0pt}%
  \let\big=\ninebig\normalbaselines\rm}

\def\eightpoint{\def\rm{\fam0\eightrm}
  \textfont0=\eightrm \scriptfont0=\sixrm
  \scriptscriptfont0\fiverm
  \textfont1=\eighti \scriptfont1=\sixi
  \scriptscriptfont1\fivei 
  \textfont2=\eightsy \scriptfont2=\sixsy
  \scriptscriptfont2\fivesy 
  \textfont3=\tenex \scriptfont3=\tenex
  \scriptscriptfont3\tenex 
  \textfont\itfam=\eightit\def\it{\fam\itfam\eightit}%
  \textfont\slfam=\eightsl\def\sl{\fam\slfam\eightsl}%
  \textfont\ttfam=\eighttt\def\tt{\fam\ttfam\eighttt}%
  \textfont\gothfam=\eightgoth\scriptfont\gothfam=\sixgoth 
  \scriptscriptfont\gothfam=\fivegoth
  \def\goth{\fam\gothfam\tengoth}
  \textfont\bffam=\eightbf\scriptfont\bffam=\sixbf
  \scriptscriptfont\bffam=\fivebf
  \def\bf{\fam\bffam\eightbf}%
  \tt\ttglue=.5em plus.25em minus.15em
  \normalbaselineskip=9pt \setbox\strutbox\hbox{\vrule
  height7pt depth2pt width0pt}%
  \let\big=\eightbig\normalbaselines\rm}

\def\twelvepoint{\def\rm{\fam0\twelverm}
  \textfont0=\twelverm\scriptfont0=\tenrm
  \scriptscriptfont0\sevenrm
  \textfont1=\twelvei\scriptfont1=\teni
  \scriptscriptfont1\seveni 
  \textfont2=\twelvesy\scriptfont2=\tensy
  \scriptscriptfont2\sevensy 
   \textfont\itfam=\twelveit\def\it{\fam\itfam\twelveit}%
  \textfont\slfam=\twelvesl\def\sl{\fam\slfam\twelvesl}%
  \textfont\ttfam=\twelvett\def\tt{\fam\ttfam\twelvett}%
  \textfont\gothfam=\twelvegoth\scriptfont\gothfam=\ninegoth 
  \scriptscriptfont\gothfam=\sevengoth
  \def\goth{\fam\gothfam\twelvegoth}
  \textfont\bffam=\twelvebf\scriptfont\bffam=\ninebf
  \scriptscriptfont\bffam=\sevenbf
  \def\bf{\fam\bffam\twelvebf}%
  \tt\ttglue=.5em plus.25em minus.15em
  \normalbaselineskip=12pt \setbox\strutbox\hbox{\vrule
  height9pt depth4pt width0pt}%
  \let\big=\twelvebig\normalbaselines\rm}

\def\tenbig#1{{\hbox{$\left#1\vbox
to8.5pt{}\right.\n@space$}}}
\def\ninebig#1{{\hbox{$\textfont0=\tenrm\textfont2=
\tensy\left#1\vbox to7.25pt{}\right.\n@space$}}}
\def\eightbig#1{{\hbox{$\textfont0=\ninerm\textfont2=
\ninesy\left#1\vbox to6.5pt{}\right.\n@space$}}}


\def\pp{{\goth p}}

\def\QQ{{\bf Q}}
\def\ZZ{{\bf Z}}
\def\FF{{\bf F}}

\def\CC{{\bf C}}

\def\zmod#1{\,\,({\rm mod}\,\,#1)}

\def\proof{\noindent{\bf Proof.}}

\def\exact#1,#2,#3{0\longrightarrow#1\longrightarrow#2
\longrightarrow#3\longrightarrow 0}

\def\today{\ifcase\month\or January \or February\or
March\or April\or May\or June\or July\or August\or
September\or October\or November\or December\fi
\space\number\day, \number\year}

\def\proclaim #1. #2\par{\medbreak
\noindent{\bf#1.\enspace}{\sl#2\par}\par\medbreak}

\magnification\magstep1\newcount\newpen\newpen=50
\def\qqlineA#1 #2 #3 {\line{\kern#1truept\vrule height0.5truept depth0.15truept 
width#2truept\hskip#3truept plus0.1truept\vrule height0.5truept depth0.05truept 
width16truept\hskip#3truept plus0.1truept\vrule height0.5truept depth0.15truept 
width#2truept\kern#1truept}}
\def\qqlineB#1 #2 #3 {\line{\kern#1truept\vrule height0.5truept depth0.15truept 
width#2truept\hskip#3truept plus0.1pt\vrule height0.5truept depth0.15truept 
width#2truept\kern#1truept}}
\def\qqlineC#1 #2 {\line{\hskip#1truept plus0.1pt\vrule height0.5truept 
depth0.15truept width#2truept\hskip#1truept plus0.1pt}}
\font\sanserifeight=cmss8 at 8truept
\font\sanserifeightfive=cmss8 scaled\magstep4
\newbox\boxtorre
\setbox\boxtorre=\vtop{\hsize=16truept\offinterlineskip\parindent=0truept
	\line{\vrule height1truept depth0.5truept width1truept\hfil\vrule 
	height1truept depth0.5truept width1truept\hfil\vrule height1truept 
	depth0.5truept width1truept\hfil\vrule height1truept depth0.5truept 
	width1truept\hfil\vrule height1truept depth0.5truept width1truept}
	\line{\vrule height1truept depth1truept width16truept}
	\vskip1.5truept
	\line{\vrule height1truept depth0truept width16truept}
	\vskip1truept
	\line{\vrule height1truept depth0truept width16truept}
	\vskip1truept
	\line{\vrule height1truept depth0truept width16truept}
	\vskip1truept
	\line{\vrule height1truept depth0truept width16truept}
	\vskip1truept
	\line{\vrule height1truept depth0truept width16truept}
	\vskip1truept
	\line{\vrule height1truept depth0truept width16truept}
	\vskip1truept
	\line{\vrule height1truept depth0truept width16truept}
	\vskip1truept
	\line{\vrule height1truept depth0truept width16truept}
	\vskip1truept
	\line{\vrule height1truept depth0truept width16truept}
	\vskip1truept
	\line{\vrule height1truept depth0truept width16truept}
	\vskip1truept}
\newbox\boxu
\setbox\boxu=\hbox{\vrule height1truept depth24truept width4truept\kern6truept
	\copy\boxtorre\kern6truept\vrule height1truept depth24truept 
	width4truept}
\newbox\boxuni
\setbox\boxuni=\vtop{\baselineskip=-1000pt\lineskip=-0.15truept
\lineskiplimit=0pt\parindent=0truept\hsize=36truept
\line{\copy\boxu\hfil}
\qqlineA 0.0 4.1 5.9 
\qqlineA 0.0 4.2 5.8 
\vskip-0.1truept
\qqlineB 0.1 4.2 27.4 
\qqlineB 0.21 4.2 27.2 
\qqlineA 0.3 4.2 5.5 
\qqlineA 0.41 4.2 5.4 
\qqlineB 0.5 4.3 26.4 
\qqlineB 0.6 4.4 26.0 
\qqlineA 0.7 4.5 4.8 
\qqlineA 0.8 4.7 4.5 
\qqlineB 1.0 4.8 24.4 
\qqlineB 1.2 4.8 24.0 
\qqlineB 1.5 4.8 23.4 
\qqlineB 1.7 4.91 22.8 
\qqlineB 1.9 5.1 22.0 
\qqlineB 2.21 5.2 21.2 
\vskip-0.1truept
\qqlineB 2.5 5.3 20.4 
\qqlineB 2.8 5.5 19.4 
\qqlineB 3.1 5.7 18.4 
\qqlineB 3.4 6.0 17.21 
\qqlineB 3.7 6.3 16.0 
\qqlineB 4.2 6.8 14.0 
\qqlineB 4.6 7.4 12.0 
\qqlineB 5.1 7.9 10.0 
\qqlineB 5.6 9.0 6.8  
\vskip-0.1truept
\qqlineC 6.1 23.8 
\qqlineC 6.7 22.6 
\qqlineC 7.4 21.21 
\qqlineC 8.0 20.0 
\qqlineC 8.6 18.8 
\vskip-0.1truept
\qqlineC 9.6 16.8
\qqlineC 10.91 14.2 
\qqlineC 12.0 12.0 
\qqlineC 13.5 9.0 
\vskip 3.1truept
\centerline{\sanserifeight TOR VERGATA}}

\newbox\boxprelim
\setbox\boxprelim\vtop to\dp\boxuni{\parindent=0truept\hsize=5truecm\eightpoint
Preliminary version\par\today\par\vfill}
\newbox\boxroma
\setbox\boxroma\vtop
to\dp\boxuni{\hsize=6truecm\parindent=0truept \vglue
.0truecm\sanserifeightfive II Universit\`a degli\par
\smallskip \hskip 1truecm Studi di Roma\par\vfill}

\def\title #1\par{
  \hbox{\copy\boxprelim\hskip 1.2in
\copy\boxuni\hskip 1cm\copy\boxroma}  
  \bigskip\bigskip\bigskip 
  \bigskip\bigskip\bigskip\bigskip\noindent
  {\twelvepoint$\!\!\!\!$\bf#1}\par
  \bigskip\bigskip\noindent
  Ren\'e Schoof
  \bigskip\eightpoint\noindent
  \vbox{
  \hbox{Dipartimento di Matematica}
  \hbox{$\hbox{2}^{\hbox{a}}$ Universit\`a di Roma 
  ``Tor Vergata"}
  \hbox{I-00133 Roma ITALY}
  \hbox{Email: \tt schoof@science.uva.nl}}
  \bigskip\bigskip
  \tenpoint
  }
\def\abstract #1\par{\eightpoint\vbox{\noindent
  {\bf Abstract.\ }#1\par}\bigskip\tenpoint
  }
\def\bibliography#1\par{\vskip0pt
  plus.3\vsize\penalty-250\vskip0pt
  plus-.3\vsize\bigskip\vskip\parskip
  \message{Bibliography}\leftline{\bf
  Bibliography}\nobreak\smallskip\noindent
  \ninepoint\frenchspacing#1}
\outer\def\beginsection#1\par{\vskip0pt
  plus.3\vsize\penalty\newpen\newpen=-50\vskip0pt
  plus-.3\vsize\bigskip\vskip\parskip
  \message{#1}\leftline{\bf#1}
  \nobreak\smallskip\noindent}

\def\AB{{1}}  
\def\ANT{{2}} 
\def\BK{{3}}  
\def\FCALE{{4}} 
\def\FAL{{6}} 
\def\F{{7}}   
\def\KM{{9}}  
\def\MA{{10}}  
\def\MAZ{{11}} 
\def\SCH{{16}}  
\def\SCHO{{17}} 
\def\OT{{19}} 
\def\SGA7{{8}}
\def\ED{{5}}   
\def\ME{{12}} 
\def\ODL{{13}}
\def\RAY{{15}}
\def\WAS{{21}}
\def\OPDV{{14}} 
\def\TDB{{18}} 
\def\TAU{{20}}

\def\gp{{G[\pi]}}

\title 
Semi-stable abelian varieties over $\QQ$ with
bad reduction in one prime only

\abstract
Let $l$ be a prime. We show that there do not exist any non-zero semi-stable abelian varieties over
$\QQ$ with good reduction outside~$l$ if and only if $l=2$, 3, 5, 7~or~13.  We show
that  any  semi-stable abelian variety over $\QQ$ with good reduction outside~$l=11$
is isogenous to a power of the Jacobian variety of the modular curve~$X_0(11)$.
In addition we show that there do not exist any non-zero abelian varieties over $\QQ$ with
good reduction outside~$l$ acquiring semi-stable reduction at~$l$ over a tamely ramified
extension if and only if $l\le 5$. 

\beginsection 1. Introduction.

In this paper we study abelian varieties over~$\QQ$ that have good reduction
at all primes except one. Denoting this one bad prime by~$l$, our first result is concerned with
abelian varieties that are semi-stable at~$l$.  Examples of these are provided by the
Jacobian varieties
$J_0(l)$ of the modular curves
$X_0(l)$. For $l=11$ or~$l\ge17$, these are  non-zero
abelian varieties. In the other direction we prove  for the 
primes~$l$ for which $J_0(l)$ is zero ---i.e. for $l=2$, $3$, $5$, $7$ and~$13$--- 
the following result.

\proclaim Theorem 1.1.
Let $l$ be a prime. There do not exist any non-zero semi-stable abelian varieties over~$\QQ$ that
have  good reduction at every prime different from~$l$ 
if and only  if $l=2$, 3, 5, 7~or~13.

\medskip\noindent
For the prime $l=11$ we show the following.

\proclaim Theorem 1.2. Every semi-stable abelian variety over~$\QQ$ that has good reduction outside
the prime~11 is isogenous to a power of $J_0(11)$.

\medskip\noindent
Jacobian varieties $J_0(l^2)$ of the modular curves $X_0(l^2)$ are
other examples of abelian varieties over~$\QQ$ that
have good reduction at all primes different from~$l$. These abelian varieties are {\it not}
semi-stable at~$l$. However, S.J.~Edixhoven~[\ED] showed that $J_0(l^2)$  acquires semi-stable
reduction at~$l$ over an extension that is merely tamely ramified at~$l$.  When $l\ge7$ the
abelian varieties $J_0(l^2)$ are not zero. In this paper we prove for the 
primes~$l$ for which $J_0(l^2)$ is zero ---i.e., for $l=2$, $3$ and~$5$--- the following theorem.

\proclaim Theorem 1.3.
Let $l$ be a prime. There do not exist any non-zero abelian varieties over~$\QQ$ that have  good
reduction at every prime different from~$l$ and that acquire semi-stable reduction at~$l$ over an
extension of~$\QQ$ that is at most tamely ramified at~$l$ if and only  if~$l=2$, $3$ or~$5$.

\medskip\noindent
For $l=7$ the dimension of $J_0(l^2)$ is~1. It is not unlikely that every abelian variety over~$\QQ$ with 
good reduction outside~$7$ and acquiring semi-stable reduction at~$7$ over an
extension that is at most tamely ramified at~7, is isogenous to a power 
of~$J_0(49)$, but I was not able to prove this.

Our results are consistent with conditional results obtained by J.-F.
Mestre~[\ME,~III]. Mestre's methods are analytic in nature. He assumes that the $L$-functions
associated to abelian varieties over~$\QQ$ admit analytic continuations to~$\CC$ and he applies
 A.~Weil's explicit formulas. Our results do not depend on any unproved hypotheses. In a
recent paper~[\BK] A.~Brumer and K.~Kramer  prove  Theorem~1.1 for primes~$l\le 7$. In
this paper we prove the somewhat stronger Theorem~1.3 for the primes $l=2$, 3 and~5,
while in addition to $l=7$,  we  take care of the primes~$11$ and~13. Our
proof, like the one by Brumer and Kramer, proceeds by studying for a suitable small
prime~$p\not= l$, the $p^n$-torsion points 
$A[p^n]$  of  abelian varieties~$A$ that have good reduction at every prime different from~$l$ and
that either have  semi-stable reduction at~$l$ or acquire it over a tamely ramified extension. When $l=2$,
3, 5, 7~or~13  we show for every
$n\ge 1$ that the group scheme $A[p^n]$  is an extension of a constant group scheme by a diagonalizable
one. This leads to a contradiction when
${\rm dim}(A)>0$. The way in which the contradiction is obtained differs  from Brumer and Kramer's
method. Our method is closer to the approach taken by J.-M. Fontaine in~[\F].

For $l=11$ things are more complicated. In this case there may exist group 
schemes of $p$-power order that
are not extensions of a constant group scheme by a diagonalizable one.  Indeed, for $p=2$, the 2-torsion
subgroup scheme of $J_0^{}(11)$ is an example of such an `exotic' group scheme.

Finally we mention F.~Calegari's result~[\FCALE]. Under assumption of the Generalized Riemann
Hypothesis, he determines all squarefree $n>0$ for which there do not exist any non-zero abelian
varieties over $\QQ$ with good reduction at all primes not dividing~$n$ while the reduction at the
primes dividing $n$ is semi-stable. These turn out to be the squarefree $n$ for which the genus
of~$X_0(n)$ is zero: $n=1,2,3,5,6,7,10$ and~13.

For any pair of distinct primes $p$ and~$l$ we introduce in section 2 two categories $\underline{C}$ and
$\underline{D}$  of finite flat group schemes of $p$-power order over the ring~$\ZZ[{1\over l}]$. In
terms of these we formulate in section~3 two simple criteria for Theorem~1.1 or~1.3 to
hold. Each criterion has two parts. One involves the  extensions of the group schemes
$\mu_p$ by~$\ZZ/p\ZZ$ over~$\ZZ[{1\over l}]$. We determine these in section~4. 
The other is concerned with the simple objects in in the categories~$\underline{D}$ and
$\underline{C}$ respectively. In sections~5 and~6 we determine these for a very short list of
pairs~$(p,l)$.  Theorems~1.1 and~1.3 follow from this. We deal with the prime~$l=11$
and prove Theorem~1.2 in section~7. Section~8 contains a theorem concerning $p$-divisible
groups that is essential for our proof of Theorem~1.2.

I thank Frank Calegari for pointing
out that the methods used to prove Theorem~1.1 could also be used to prove Theorem~1.3.

\beginsection 2. Two categories of finite flat group schemes over~$\ZZ[{1\over l}]$. 

Let $p$ and $l$ be two distinct primes. By $\underline{Gr}$
we denote the category of finite flat commutative $p$-power order group schemes, or {\it
$p$-group schemes} for short,  over~$\ZZ[{1\over l}]$.
In this section we introduce two full subcategories of $\underline{Gr}$. In terms of these we
formulate in the next section two criteria for Theorems~1.1 and~1.3 to hold. 

\medskip\noindent{\bf Definition 2.1.}\ Let $p$ and $l$ be two distinct primes.
\smallskip\noindent {\sl (i)}\  The category $\underline{C}$ is the
full subcategory of $\underline{Gr}$ of group schemes $G$ for which the inertia group of every prime
over~$l$ acts tamely on~$G(\overline{\QQ})$. 
\smallskip\noindent {\sl (ii)}\ 
The category $\underline{D}$ is  the full subcategory of $\underline{Gr}$ of group schemes $G$ for which
we have  that $(\sigma-{\rm id})^2=0$ on~$G(\overline{\QQ})$ for all $\sigma$ in  the
inertia group $I_{\goth l}$ of any of the primes $\goth l$ lying over~$l$.

\medskip\noindent 
Since every group scheme $G$ in the category~$\underline{D}$ has $p$-power order, the  relation
$(\sigma-{\rm id})^2=0$ implies that
$\sigma^{p^n}={\rm id}$ for some $n\ge0$. This implies that every inertia
group $I_{\goth l}$ acts through a finite $p$-group on~$G(\overline{\QQ})$.  Since $p\not=l$, this action
is tame and hence  
$\underline{D}$ is a full subcategory of~$\underline{C}$:
$$
\underline{D}\quad\subset\quad\underline{C}\quad\subset\quad\underline{Gr}.
$$
We now exhibit several objects in the categories $\underline{C}$ and $\underline{D}$. The first
examples explain the very definition of these two categories.
\smallskip
\noindent {\bf 1.\ } By \hyphenation{Gro-then-dieck} A.~Grothendieck's semi-stable
reduction Theorem~[\SGA7, Exp.IX,~(3.5.3)], the subgroup schemes
$A[p^n]$ of $p^n$-torsion points of semi-stable abelian varieties $A$ over~$\QQ$ that have good
reduction outside~$l$ are objects of~$\underline{D}$  and hence of~$\underline{C}$. On the other
hand, the subgroup schemes $A[p^n]$ of abelian varieties $A$ over~$\QQ$ that have good reduction
outside~$l$ and acquire semi-stable reduction at~$l$ over some tamely ramified extension, are objects
of~$\underline{C}$ and need not be objects of~$\underline{D}$. 

\smallskip\noindent {\bf 2.\ } 
Constant and diagonalizable group schemes of $p$-power order are objects
of~$\underline{D}$ and hence of~$\underline{C}$.  So are certain twists of these group schemes that
are unramified outside~$l$. Here are some explicit examples: let $n\ge 0$ and $V=(\ZZ/p\ZZ)^n$. For any
representation
$\varrho:{\rm Gal}(\overline{\QQ}/\QQ)\longrightarrow{\rm GL}(V)$ we obtain a Galois module with
underlying group $V$. If $\varrho$ is unramified outside $l$ and infinity, this module is
the group of points of a group scheme $V(\varrho)$ that is \'etale over~$\ZZ[{1\over l}]$. 
If
$\rho$ at most tamely ramified at~$l$, the group scheme $V(\varrho)$ is an object
of~$\underline{C}$. 

For intance, taking $V=(\ZZ/p\ZZ)^{l-1}$ and $\varrho:{\rm Gal}(\QQ(\zeta^{}_l)/\QQ)\longrightarrow{\rm
GL}(V)$ the permutation representation, we obtain an object of~$\underline{C}$ whose points generate the
field $\QQ(\zeta_l^{})$. Another example is given by $V=(\ZZ/p\ZZ)^2$ and $\varrho$ a homomorphism 
$$
\varrho:{\rm Gal}(\QQ(\zeta_l^{})/\QQ)\quad\longrightarrow\quad
\{\hbox{$\pmatrix{1&x\cr0&1\cr}$}:x\in\FF_p\}\quad\subset\quad{\rm GL}(V).
$$
This time the  group scheme $V(\varrho)$ is an object in~$\underline{D}$, because $(\varrho(\sigma)-{\rm
id})^2=0$ for all~$\sigma$. For this twist to be non-trivial, it is necessary that $l\equiv 1\zmod p$.
Taking Cartier duals we obtain unramified twists of diagonalizable group schemes that are objects
in $\underline{C}$ and~$\underline{D}$ respectively.

\smallskip\noindent {\bf 3.\ }  Cartier duals of objects in~$\underline{D}$ are also
in~$\underline{D}$. Any closed flat  subgroup scheme and any quotient by such a group scheme of an 
object in $\underline{D}$ is again in
$\underline{D}$. The product of any two objects in $\underline{D}$ is again in~$\underline{D}$.
It follows from the definition of the Baer sum that if $G_1$ and $G_2$  are objects in $\underline{D}$,
then  the extension classes of
$G_1$ by
$G_2$ that are themselves objects of~$\underline{D}$  make up a {\it subgroup}
$$
{\rm Ext}^1_{\underline{D}}(G_1,G_2)
$$ 
of the group ${\rm Ext}^1_{\ZZ[{1\over l}]}(G_1,G_2)$  of {\it all} extensions of $G_1$
by~$G_2$  in the category~$\underline{Gr}$.
All these remarks also hold for the category~$\underline{C}$, but
an extension $G$ of an object $G_1\in\underline{C}$ by another
$G_2\in\underline{C}$ is {\it automatically} an object of $\underline{C}$. In other words ${\rm
Ext}^1_{\underline{C}}(G_1,G_2)={\rm Ext}_{\ZZ[{1\over l}]}^1(G_1,G_2)$. Indeed, there is
some  exponent $e$ that is prime to~$l$ and has the property that  for all primes $\goth l$
over~$l$ and all 
$\sigma\in I_{\goth l}$, the automorphism $\sigma^e$ acts
trivially on the points of $G_1$ and~$G_2$. It follows that for some $s\ge 0$, the
automorphism
$\sigma^{ep^s}$ acts trivially on the points of~$G$. This implies that any inertia group $I_{\goth l}$
acts through its tame quotient. 

It is in general not true that ${\rm
Ext}^1_{\underline{D}}(G_1,G_2)={\rm Ext}^1_{\ZZ[{1\over l}]}(G_1,G_2)$.
It is true however when the inertia groups $I_{\goth l}$ act trivially
on the points of $G_1$ and~$G_2$. We have for instance that 
$${\rm Ext}^1_{\underline{D}}(\mu_p,\ZZ/p\ZZ)={\rm
Ext}^1_{\underline{C}}(\mu_p,\ZZ/p\ZZ)={\rm Ext}^1_{\ZZ[{1\over l}]}(\mu_p,\ZZ/p\ZZ).
$$
Since $\underline{D}$
is closed under the formation of products, closed subgroup schemes and quotients by closed
flat subgroup schemes, there is for every short exact sequence $0\rightarrow G_1\rightarrow
G_2\rightarrow G_3\rightarrow 0$ in~$\underline{D}$ and any $H\in\underline{D}$ a six term 
exact sequence
$$\eqalign{
0\longrightarrow\,&{\rm Hom}(H,G_1)\longrightarrow{\rm Hom}(H,G_2)\longrightarrow{\rm
Hom}(H,G_3)\longrightarrow\cr
\longrightarrow\,&{\rm Ext}^1_{\underline{D}}(H,G_1)\longrightarrow{\rm
Ext}^1_{\underline{D}}(H,G_2)\longrightarrow{\rm Ext}^1_{\underline{D}}(H,G_3)\cr}
$$
The contravariant version of this statement is also true.

\smallskip\noindent{\bf 4.\ }
N.~Katz and B.~Mazur construct in their book~[\KM, Interlude~8.7] some explicit extensions of
$\ZZ/p\ZZ$ by~$\mu_p$.  These are objects in~$\underline{D}$. We specialize their
construction to our situation. For any unit
$\varepsilon\in\ZZ[{1\over l}]^*$ they define a finite flat $p$-group scheme $G_{\varepsilon}$
over
$\ZZ[{1\over l}]$ of order~$p^2$.  It is killed by~$p$ and fits in
an exact sequence
$$
0\longrightarrow\mu_p\longrightarrow G_{\varepsilon}\longrightarrow\ZZ/p\ZZ\longrightarrow 0.
$$
Two group schemes $G_{\varepsilon}$ and
$G_{\varepsilon'}$ are isomorphic if and only if
$\varepsilon/\varepsilon'=u^p$ for some $u\in\ZZ[{1\over l}]^*$.
The points of $G_{\varepsilon}$ generate the number field $\QQ(\zeta_p,{\root p\of{\varepsilon}})$.
For $p=2$ and $\varepsilon=-1$, we recover the group scheme $D$ in~[\MAZ,~Prop.4.2]. 

\medskip We conclude this section by giving an explicit example of an object
in~$\underline{D}$ and an extension of it by itself that is {\it not}
in~$\underline{D}$. It plays no  role in the rest of this paper. See section~7 for a
similar example. We take $l=7$ and $p=2$. The
elliptic curve $E$ given by the equation $Y^2+XY=X^3-X^2-2X-1$ has conductor~49.
See~[\ANT].  Its endomorphism ring over~$\QQ(\sqrt{-7})$ is isomorphic to the
ring~$\ZZ[{{1+\sqrt{-7}}\over2}]$. The points
$\infty$ and $(2,-1)$  are rational while the other two 2-torsion points are
only defined over~$\QQ(\sqrt{-7})$. Therefore 
the non-trivial automorphism of~$\QQ(\sqrt{-7})$ acts as the matrix
$A=\pmatrix{1&1\cr0&1\cr}$ on the  points of the group scheme~$E[2]$ of 2-torsion points. Since $A-{\rm
id}$ has square zero, the group scheme
$E[2]$ is an object in~$\underline{D}$.  The group scheme $E[4]$ is an extension of $E[2]$
by~$E[2]$. The points of~$E[4]$ generate the field $L=\QQ(i,{\root 4\of {-7}})$.
There is only one prime over~7 in $L$ and its inertia group is ${\rm Gal}(L/\QQ(i))$. It is
cyclic of order~4. The square of a generator $\sigma$ commutes with the endomorphisms of~$E$.
Therefore it respects the eigenspaces of the squares of the two endomorphism of~$E$ of degree~2.
Since its determinant is a square in $(\ZZ/4\ZZ)^*$, it is equal to~$1$. It follows that
$\sigma^2$  acts as a scalar matrix on $E[4]$. If it were true that
$(\sigma-{\rm id})^2=0$, then $2\sigma=\sigma^2+{\rm id}$ on~$E[4]$ and hence $\sigma$ would act a
scalar on $E[2]$. Since this is not the case, the group scheme $E[4]$ is not an object
in the category~$\underline{D}$.

\beginsection 3. Two criteria.

In this section we formulate criteria for Theorems~1.1 and~1.3 to hold. Let $l$ and $p$
be two distinct primes. A {\it simple} object in any of the categories $\underline{Gr}$,
$\underline{C}$ or
$\underline{D}$ of the previous section is an object in that category
that does not admit any non-trivial closed flat
subgroup schemes.

\proclaim Proposition 3.1.  Let $l$ be a prime and suppose there exists a prime $p\not=l$ for
which 
\smallskip
\item{--}  the only simple objects in the category $\underline{D}$ are $\ZZ/p\ZZ$ and~$\mu_p$;
\smallskip
\item{--}  ${\rm Ext}^1_{\ZZ[{1\over l}]}(\mu_p,\ZZ/p\ZZ)=0$.
\smallskip\noindent Then
there do not exist any non-zero abelian varieties over $\QQ$ that have good reduction at every prime
different from~$l$ and are semi-stable at~$l$.

\smallskip
\proof\  The proof is very similar to the one in~[\F,~section~3.4.3]. We briefly sketch it. Let $A$ be a
$g$-dimensional  abelian variety  over
$\QQ$ that has  good reduction at every prime different from~$l$ and semi-stable reduction at~$l$. Let
$n\ge 1$ and consider the  closed subgroup scheme $A[p^n]$ of $p^n$-torsion points
over~${\ZZ[\hbox{${1\over l}$}]}$. This is an object of the category~$\underline{D}$. We
filter
$A[p^n]$ with flat closed subgroup schemes and successive {\it simple} subquotients. The simple
steps are by assumption isomorphic to $\ZZ/p\ZZ$ and~$\mu_p$. By the second assumption, any extension
$$
0\quad\longrightarrow\quad\ZZ/p\ZZ\quad\longrightarrow\quad G
\quad\longrightarrow\quad\mu_p \quad\longrightarrow\quad 0
$$
splits over~$\ZZ[{1\over l}]$. Therefore we can modify the filtration and obtain
for every $n\ge 1$ an exact sequence of finite flat $\ZZ[{1\over l}]$-group schemes
$$
0\quad\longrightarrow\quad M_n\quad\longrightarrow \quad A[p^n]\quad
\longrightarrow \quad C_n\quad\longrightarrow \quad 0 
$$
with $M_n$ an extension of group schemes isomorphic to $\mu_p$ and 
$C_n$ an extension of copies of~$\ZZ/p\ZZ$.
The fundamental group of~$\ZZ[{1\over l}]$ acts on the points of~$C_n$, and it does so through a finite
$p$-group~$P$. Since the maximal abelian extension of~$\QQ$ that is unramified
outside~$l\cdot\infty$ is contained in the cyclic extension $\QQ(\zeta_l^{})$,  the group
$P/P'$ is cyclic. This implies that $P$ itself is cyclic. It follows that $C_n$ becomes {\it
constant} over the ring~$\ZZ[{1\over l},\zeta_l^{}]$. Similarly, it follows from Cartier
duality that $M_n$ becomes diagonalizable over this ring.

Pick a non-zero prime $\pp$ of~$\ZZ[{1\over l},\zeta_l^{}]$ and let $k_{\pp}$ denote its
residue field. The abelian variety
$A/M_n$ has at least $\#C_n$ rational points modulo~$\pp$. Since $A/M_n$ is isogenous
to~$A$, it has the same number of points as~$A$ modulo~$\pp$. This implies that
$\#C_n\le\#A(k_{\pp})$. Similarly, the abelian variety $A^{\rm
dual}/C_n^{\vee}$ has at least
$\#M_n^{\vee}$ points modulo~$\pp$. Here $G^{\vee}$ denotes the Cartier dual of a finite
group scheme~$G$. Since  $A^{\rm dual}/C_n^{\vee}$ is also isogenous to $A$, we see that 
$\#M_n=\#M_n^{\vee}\le\#A(k_{\pp})$. However, the product of the orders of $M_n$ and
$C_n$ is equal to $\#A[p^n]=p^{2ng}$. This leads
to a contradiction when $n\rightarrow\infty$ unless
$g=0$. This proves the proposition.

\proclaim Proposition 3.2.   Let $l$ be a prime and suppose there exists a prime $p\not=l$
for which 
\smallskip
\item{--}  the only simple objects in the category $\underline{C}$ are $\ZZ/p\ZZ$ and~$\mu_p$;
\smallskip
\item{--}  ${\rm Ext}^1_{\ZZ[{1\over l}]}(\mu_p,\ZZ/p\ZZ)=0$.
\smallskip\noindent Then there do not exist
any non-zero abelian varieties over $\QQ$ that have good reduction at every prime different from~$l$ and
acquire semi-stable reduction at~$l$ over a tamely ramified extension.

\proof\ The proof is similar to the proof of the previous proposition. This time the
group scheme~$A[p^n]$ is an object of the category~$\underline{C}$. We leave the proof to the reader.

\beginsection 4. Extensions of $\mu_p$ by $\ZZ/p\ZZ$ over~$\ZZ[{1\over l}]$.

The criteria of the previous section involve the group
${\rm Ext}^1(\mu_p,\ZZ/p\ZZ)$ and the simple objects in the categories
$\underline{D}$ or~$\underline{C}$.
We discuss the simple objects in the next section. In this section we fix two distinct
primes $l$ and~$p$ and  determine the extensions of the group scheme $\mu_p$ by $\ZZ/p\ZZ$ 
over the ring~$\ZZ[{1\over l}]$. For $p=2$ see
also~[\MAZ,~Prop.5.1].

\proclaim Proposition 4.1.  Let $l$ and $p$ be distinct primes. We have that
$$
{\rm dim}^{}_{\FF_p}{\rm Ext}^1_{\ZZ[{1\over l}]}(\mu_p,\ZZ/p\ZZ)=\cases{1,&if
${{l^2-1}\over{24}}\equiv 0\zmod p$;\cr0,&otherwise.\cr}
$$

\proof\ Note  that ${{l^2-1}\over{24}}$ is $p$-integral since $p\not=l$.
The condition ${{l^2-1}\over{24}}\equiv 0\zmod{p}$ is
just a compact way of saying that $l\equiv\pm 1\zmod{p}$ when $p\ge 5$, that $l\equiv\pm
1\zmod{9}$ when $p=3$ and  that $l\equiv\pm 1\zmod{8}$ when $p=2$. We work with the following base rings
and invoke the Mayer-Vietoris exact sequence of~[\SCH,~Prop.2.4]. 
$$
\def\normalbaselines{\baselineskip20pt\lineskip3pt 
\lineskiplimit3pt}
\matrix{&&\ZZ_p&&\cr&\nearrow&&\searrow&\cr
\ZZ[\hbox{${1\over l}$}]&&&&\QQ_p\cr
&\searrow&&\nearrow&\cr&&\ZZ[\hbox{${1\over{lp}}$}]&&\cr}
$$
Since $\mu_p$ is connected, while $\ZZ/p\ZZ$ is \'etale, the group ${\rm
Hom}_{\ZZ_p}(\mu_p,\ZZ/p\ZZ)$ vanishes. Therefore ${\rm Hom}_{\ZZ[{1\over
l}]}(\mu_p,\ZZ/p\ZZ)$ is zero as well. The natural homomorphism ${\rm
Hom}_{\ZZ_p[{1\over{lp}}]}(\mu_p,\ZZ/p\ZZ)\longrightarrow{\rm Hom}_{\QQ_p}(\mu_p,\ZZ/p\ZZ)$
is an isomorphism. More precisely, both groups are zero when $p$ is odd, while they both
have order~2 when~$p=2$. Finally, the group ${\rm Ext}^1_{\ZZ_p}(\mu_p,\ZZ/p\ZZ)$ is trivial
because any extension of $\mu_p$ by $\ZZ/p\ZZ$ over~$\ZZ_p$ is split by the connected
component. Therefore the Mayer-Vietoris sequence provides us with
the following exact sequence. Since the group schemes $\mu_p$ 
and $\ZZ/p\ZZ$ are \'etale over the rings 
$\ZZ[{1\over{pl}}]$ and $\QQ_p$, this sequence gives a description of the
group ${\rm Ext}^1_{\ZZ[{1\over l}]}(\mu_p,\ZZ/p\ZZ)$ in terms of groups of
extensions of Galois modules:
$$
0\longrightarrow{\rm Ext}^1_{\ZZ[{1\over l}]}(\mu_p,\ZZ/p\ZZ)\longrightarrow {\rm
Ext}^1_{\ZZ[{1\over{pl}}]}(\mu_p,\ZZ/p\ZZ)\longrightarrow {\rm
Ext}^1_{\QQ_p}(\mu_p,\ZZ/p\ZZ). 
$$
Next we adjoin the $p$-th roots of unity. Let $\Delta$ denote the Galois group of
$\ZZ[{1\over{pl}},\zeta_p]$ over~$\ZZ[{1\over{pl}}]$.  It acts naturally on the group 
${\rm Ext}^1_{\ZZ[{1\over{pl}},\zeta_p]}(\ZZ/p\ZZ,\mu_p)$ of extensions of $\mu_p$ by~$\ZZ/p\ZZ$
over~$\ZZ[{1\over{pl}},\zeta_p]$. Since
$\Delta$ has order prime to~$p$, the extensions  over the ring
$\ZZ[{1\over{pl}}]$ correspond precisely to the ones that are $\Delta$-invariant. Since
$\ZZ[{1\over{pl}},\zeta_p]$ contains $\zeta_p$, the groups ${\rm
Ext}^1_{\ZZ[{1\over{pl}},\zeta_p]}(\mu_p,\ZZ/p\ZZ)$ and ${\rm
Ext}_{\ZZ[{1\over{pl}},\zeta_p]}^1(\ZZ/p\ZZ,\mu_p)$ are isomorphic! However, 
under this isomorphism the
$\Delta$-invariant extensions of
$\mu_p$ by~$\ZZ/p\ZZ$ correspond to the extensions that are contained in the
$\omega^2$-eigenspace ${\rm Ext}^1_{\ZZ[{1\over{pl}},\zeta_p]}(\ZZ/p\ZZ,\mu_p)_{\omega^2}$.  
Mutatis mutandis, the same holds for $\QQ_p$ and
its extension $\QQ_p(\zeta_p)$. Here
$\omega:\Delta\longrightarrow\ZZ^*_p$ denotes the 
cyclotomic (or Teichm\"uller) character and
$M_{\omega^2}$ denotes the
$\omega^2$-eigenspace of a $\ZZ_p[\Delta]$-module~$M$. 
These considerations lead to the exact sequence
$$
0\longrightarrow{\rm Ext}^1_{\ZZ[{1\over l}]}(\mu_p,\ZZ/p\ZZ)\longrightarrow {\rm
Ext}^1_{\ZZ[{1\over{pl}},\zeta_p]}(\ZZ/p\ZZ,\mu_p)_{\omega^2}\longrightarrow {\rm
Ext}^1_{\QQ_p(\zeta_p)}(\ZZ/p\ZZ,\mu_p)_{\omega^2}.
$$
We use the exact sequence
$0\rightarrow\ZZ\rightarrow\ZZ\rightarrow\ZZ/p\ZZ\rightarrow 0$ to express the two groups ${\rm
Ext}^1_{\ZZ[{1\over{pl}},\zeta_p]}(\ZZ/p\ZZ,\mu_p)$ and ${\rm
Ext}^1_{\QQ_p(\zeta_p)}(\ZZ/p\ZZ,\mu_p)$
in terms of Galois cohomology. This leads to the following 
commutative diagram of
$\FF_p[\Delta]$-modules.
$$
\def\normalbaselines{\baselineskip20pt\lineskip3pt 
\lineskiplimit3pt}
\matrix{0&\longrightarrow&\mu_p&\longrightarrow {\rm
Ext}^1_{\ZZ[{1\over{lp}},\zeta_p]}(\ZZ/p\ZZ,\mu_p)&\longrightarrow&
H^1(G_{\ZZ[\hbox{$1\over{lp}$},\zeta_p]},\mu_p)
&\longrightarrow 0\cr
&&\Vert&\Big\downarrow&& \Big\downarrow&&\cr
0&\longrightarrow&\mu_p&\longrightarrow {\rm
Ext}^1_{\QQ_p(\zeta_p)}(\ZZ/p\ZZ,\mu_p)&\longrightarrow&
H^1(G_{\QQ_p(\zeta_p)},\mu_p)&\longrightarrow 0\cr}
$$
Here $G_{\QQ_p(\zeta_p)}$ denotes the absolute Galois group of
$\QQ_p(\zeta_p)$ and $G_{\ZZ[\hbox{$1\over{lp}$},\zeta_p]}$ is the fundamental group of
$\ZZ[\hbox{$1\over{lp}$},\zeta_p]$. It follows that there is an exact sequence of
$\FF_p[\Delta]$-modules.
$$
0\longrightarrow{\rm Ext}^1_{\ZZ[{1\over l}]}(\mu_p,\ZZ/p\ZZ)\longrightarrow 
H^1(G_{\ZZ[\hbox{$1\over{lp}$},\zeta_p]},\mu_p)_{\omega^2}\longrightarrow
H^1(G_{\QQ_p(\zeta_p)},\mu_p)_{\omega^2}.
$$
The Kummer sequences over the rings $\ZZ[\hbox{$1\over{lp}$},\zeta_p]$ and $\QQ_p(\zeta_p)$
give rise to the following commutative diagram of
$\FF_p[\Delta]$-modules
$$
\def\normalbaselines{\baselineskip20pt\lineskip3pt 
\lineskiplimit3pt}
\matrix{0\longrightarrow&\!\!\ZZ[\hbox{$1\over{lp}$},\zeta_p]^*/(\ZZ[\hbox{$1\over{lp}$},\zeta_p]^*)^p
&\!\longrightarrow&
\!H^1(G_{\ZZ[\hbox{$1\over{lp}$},\zeta_p]},\mu_p)&\! \longrightarrow&
\!\! Cl(\ZZ[\hbox{$1\over{lp}$},\zeta_p])[p]
&\!\longrightarrow 0\cr
&\!\Big\downarrow&&\!\Big\downarrow&&\cr
&\!\QQ_p(\zeta_p)^*/(\QQ_p(\zeta_p)^*)^p&\!\mathop{\longrightarrow}\limits^{\cong} &\!
H^1(G_{\QQ_p(\zeta_p)},\mu_p)&&\cr}
$$
We take $\omega^2$-eigenspaces. By Herbrand's Theorem~[\WAS,~Thm.6.17] the
$\omega^2$-eigenspace of the
$p$-part of the class group of the ring~$\ZZ[\zeta_p]$ vanishes. This implies that the
$\omega^2$-eigenspace of 
$Cl(\ZZ[\hbox{$1\over{lp}$},\zeta_p])$ is also trivial. An application of the Snake Lemma  leads then to
the following exact sequence, suitable for explicit calculation.
$$
(*)\,\,\,0\longrightarrow{\rm Ext}^1_{\ZZ[{1\over l}]}(\mu_p,\ZZ/p\ZZ)\longrightarrow
\left(\ZZ[\hbox{$1\over{lp}$},\zeta_p]^*/(\ZZ[\hbox{$1\over{lp}$},
\zeta_p]^*)^p\right)_{\omega^2}\longrightarrow
\left(\QQ^{}_p(\zeta^{}_p)^*/(\QQ^{}_p(\zeta^{}_p)^*)^p\right)^{}_{\omega^2}
$$
When $p=2$, the group $\Delta$ is trivial and $\omega=1$. The middle and rightmost 
groups in the exact sequence~($*$) are each of dimension~3 over~$\FF_2$. The middle one is generated by 
$2$, $-1$ and $l$ while the rightmost group is generated by~$2$, $-1$ and $5$. This implies that the
extension group on the left is cyclic of order~2 when $\pm l$ is a 2-adic square, i.e. when $l\equiv \pm
1\zmod 8$, while it is trivial otherwise.

When $p=3$, the group $\Delta$ has order~2, so that the $\omega^2$-eigenspace is simply the group of
$\Delta$-invariants. The middle and  rightmost  groups in the sequence~($*$)
are each of dimension~2 over~$\FF_3$. The middle one is generated by~3 and~$l$, while the rightmost one
is generated by~3 and the unit~4. This implies that the  group ${\rm Ext}^1_{\ZZ[{1\over
l}]}(\mu_3,\ZZ/3\ZZ)$ is a 1-dimensional $\FF_3$-vector space when
$l$ is a 3-adic cube, i.e. when $l\equiv\pm 1\zmod 9$, while it is zero otherwise.

When $p\ge 5$, we first compute the group in the middle of the exact sequence~($*$). Consider the natural
exact sequence
$$
0\longrightarrow \ZZ[\hbox{$1\over{p}$},\zeta_p]^*\longrightarrow
\ZZ[\hbox{$1\over{lp}$},\zeta_p]^*\mathop{\longrightarrow}\limits^{v}\mathop{\oplus}\limits_{{\goth
l}|l}\ZZ\longrightarrow Cl(\ZZ[\hbox{$1\over{p}$},\zeta_p])\longrightarrow
Cl(\ZZ[\hbox{$1\over{lp}$},\zeta_p])\longrightarrow 0.
$$
Here $v$ is the map that sends a unit $\varepsilon\in \ZZ[\hbox{$1\over{lp}$},\zeta_p]^*$ to its
valuations at the primes $\goth l$ of $\ZZ[\zeta_p]$ that lie over~$l$.
We
tensor with $\ZZ_p$ and take $\omega^2$-eigenspaces. By Herbrand's Theorem, the
$\omega^2$-eigenspace of the
$p$-part of the class group of~$\ZZ[\hbox{$1\over{p}$},\zeta_p]$ is trivial. Therefore we obtain a three
term exact sequence. It is $\ZZ_p$-split, because the rightmost term is free over~$\ZZ_p$. Taking
quotients by
$p$-th powers we obtain therefore the following exact sequence of $\omega^2$-eigenspaces.
$$
0\longrightarrow \left(\ZZ[\hbox{$1\over{p}$},\zeta_p]^*/(\ZZ[\hbox{$1\over{p}$},
\zeta_p]^*)^p\right)_{\omega^2}\longrightarrow
\left(\ZZ[\hbox{$1\over{lp}$},\zeta_p]^*/(\ZZ[\hbox{$1\over{lp}$},
\zeta_p]^*)^p\right)_{\omega^2}\mathop{\longrightarrow}\limits^{v}
\left(\mathop{\oplus}\limits_{{\goth
l}|l}\FF_p\right)_{\omega^2}\longrightarrow 0.
$$
Identifying the Galois group $\Delta$ with $(\ZZ/p\ZZ)^*$ via its action on $\mu_p$, 
the $\FF_p[\Delta]$-module $\ZZ[\hbox{$1\over{p}$},\zeta_p]^*/(\ZZ[\hbox{$1\over{p}$},
\zeta_p]^*)^p$ is isomorphic to $\mu_p\times\FF_p[\Delta/\langle-
1\rangle]$ by~[\WAS,~Prop.8.13]. On the other hand, the module
${\oplus}_{{\goth l}|l}\FF_p$ is isomorphic to $\FF_p[\Delta/\langle l\rangle]$. This shows that  the
group in the middle has dimension 1 or 2 over~$\FF_p$.
Since $p\ge 5$,  the $\omega^{2}$-eigenspace of 
$\QQ_p(\zeta_p)^*/(\QQ_p(\zeta_p)^*)^p$ has dimension~1. This follows from a short
computation. By~[\WAS,~Thm.8.25], the $\omega^2$-eigenspace of the cyclotomic units in
$\ZZ[\hbox{$1\over{p}$},\zeta_p]^*$ maps surjectively onto it.
It follows that the rightmost arrow in  the exact sequence ($*$) is surjective and that
${\rm Ext}^1_{\ZZ[{1\over l}]}(\mu_p,\ZZ/p\ZZ)$ has dimension 1 or 0 over~$\FF_p$ depending
on whether the $\omega^{2}$-eigenspace of
${\oplus}_{{\goth
l}|l}\FF_p$ is non-trivial or not. Since $\omega$ generates the group of characters
$\Delta\longrightarrow \FF_p^*$, this depends on whether
$l\equiv
\pm1\zmod p$ or not.

This proves the proposition.

\beginsection 5. Simple objects in the categories $\underline{C}$ and $\underline{D}$.

Let $l$ and $p$ be distinct primes. In this section we give two criteria for $\ZZ/p\ZZ$ and $\mu_p$ to
be the only {\it simple} objects in the categories~$\underline{C}$ and $\underline{D}$ of
$p$-group schemes over $\ZZ[\hbox{${1\over l}$}]$ that were introduced in section~2. For any
prime~$q$, the  $q$-adic valuation
$v_q$ is normalized by putting $v_q(q)=1$.

\proclaim Proposition 5.1. Let $l$ and $p$ be distinct primes. 
When $p$ divides $l-1$, we let $F$ denote the degree $p$ subfield $F$ of~$\QQ(\zeta_l^{})$.
When not, we put $F=\QQ$. Suppose that there does not exist a Galois extension~$L$
of~$\QQ$ for which each of the following four conditions holds.
\smallskip
\item{--} the field $F(\zeta_{2p},{\root p\of l})$ is contained in~$L$;\smallskip
\item{--} the extension $F(\zeta_{2p},{\root p\of l})\subset L$ is unramified at all
primes not lying over~$p$;\smallskip
\item{--}  the $p$-adic valuation of  the root discriminant
$\delta_L$ of~$L$ is strictly smaller than~$1+{1\over{p-1}}$;\smallskip
\item{--} the degree $[L:\QQ(\zeta_p^{})]$ is not a $p$-power;\smallskip
\noindent
Then the only {\it simple} objects in the
category~$\underline{D}$ are the group schemes $\ZZ/p\ZZ$ and~$\mu_p$.

\proof\  Let $G$ be a simple object of~$\underline{D}$. When $p=2$, we put $G'=G\times
G_l\times G_{-1}$. Here $G_l$ and $G_{-1}$ denote the 
Katz-Mazur group schemes of section~2.   When $p\not=2$ we put $G'=G\times G_l\times
V(\varrho)$.  Here  $V(\varrho)$ is the twisted constant group scheme introduced in section~2 with
$V=(\ZZ/p\ZZ)^2$ and 
$$
\varrho:{\rm Gal}(\QQ(\zeta_l^{})/\QQ)\longrightarrow\{\pmatrix{1&x\cr0&1\cr}:x\in\FF_p\}
\quad\subset\quad{\rm GL}(V)
$$
of order $p$ when $l\equiv1\zmod p$ and trivial otherwise. The points of~$V(\varrho)$ generate the
field~$F$. Let $L$ be the extension of $\QQ$ generated
by the points of~$G'$. It is Galois over~$\QQ$ and by construction the first condition of the
proposition is satisfied.

Since the group scheme $G'$ is a product of objects in~$\underline{D}$, it is itself an object
in~$\underline{D}$. Therefore we know that
$(\sigma-{\rm id})^2=0$ on $G'(\overline{\QQ})$ for every
$\sigma$ in the inertia subgroup $I_{\goth l}$ of any of the primes ${\goth l}$ over~$l$. 
Since $G$ is simple, it is killed by~$p$. Therefore $G'$ is also killed by~$p$ and we have that
$\sigma^p={\rm id}$
on~$G'(\overline{\QQ})$. Since tame ramification groups are cyclic, it follows that the
ramification indices of  the primes over~$l$ divide~$p$. Since the ramification
indices of the primes over~$l$ in the subfield
$F(\zeta_{2p},{\root p\of l})$ are actually {\it equal} to~$p$, 
the field $L$ must be unramified over $F(\zeta_{2p},{\root p\of l})$ and hence the second
condition is satisfied. The estimates of
Abra\v skin and Fontaine~[\AB,~\F] for the ramification at the prime~$p$ imply that the third
condition is satisfied. Therefore $[L:\QQ(\zeta_p)]$ is a power of~$p$. 

The subgroup ${\rm Gal}(\overline{\QQ}/\QQ(\zeta_p))$ of ${\rm
Gal}(\overline{\QQ}/\QQ)$ acts on the group of points
$G(\overline{\QQ})$ through the finite $p$-group~${\rm Gal}(L/\QQ(\zeta^{}_p))$. Since
$G(\overline{\QQ})$ is a simple Galois module of
$p$-power order, it  is fixed by~${\rm Gal}(L/\QQ(\zeta^{}_p))$.
Therefore ${\rm Gal}(\overline{\QQ}/\QQ)$  acts on $G(\overline{\QQ})$ through the group 
$\Delta={\rm Gal}(\QQ(\zeta_p)/\QQ)$ of order~$p-1$.  The group $G(\overline{\QQ})$
is a product of eigenspaces. Since $G$ is simple and since  the
$p-1$-st roots of unity are in $\FF_p$, the group $G(\overline{\QQ})$ is itself equal to one
of the eigenspaces and has dimension 1 over~$\FF_p$. It follows that
$G$ has order~$p$. Since $p$ is prime in the ring $\ZZ[{1\over l}]$, the classification of
Oort-Tate~[\OT]  implies that
$G\cong\ZZ/p\ZZ$ or $\mu_p$ possibly twisted by a character $\chi$ that is
unramified outside~$l$. Such a character necessarily has order dividing~$p-1$. 
However since $G$ is an object of~$\underline{D}$, the ramification index
at~$l$ must be a power of~$p$ as well. Therefore $\chi$ is an everywhere unramified
character of~$\QQ$ and is hence trivial. This proves the theorem.

\medskip
The following theorem is a variant of Proposition~5.1. Although the conditions appear to be
similar to those of Prop.~5.1, they are much stronger. In the course of the proof we'll see
that they imply that either $l$ or~$p$ is equal to~2. In section~6 we apply 
Proposition~5.2 only to the pairs
$(l,p)=(2,3)$, $(3,2)$ and~$(5,2)$.

\proclaim Proposition 5.2. 
Let $l$ and $p$ be distinct primes. 
Suppose that there does not exist a Galois extension~$L$
of~$\QQ$ for which the following five conditions hold.
\smallskip
\item{--}  the field $\QQ(\zeta_{2p},\zeta_l^{},{\root p\of l})$ is contained in~$L$;\smallskip
\item{--}   the field $L$ is unramified outside $pl$ and infinity.  \smallskip
\item{--}  the $l$-adic valuation of  the root discriminant
$\delta_L$ of~$L$ is strictly smaller than~$1$;\smallskip
\item{--}  the $p$-adic valuation of  
$\delta_L$ of~$L$ is strictly smaller than~$1+{1\over{p-1}}$;\smallskip
\item{--} the degree $[L:\QQ(\zeta_p^{})]$ is not a $p$-power;\smallskip
\noindent
Then the only {\it simple} objects in the
category~$\underline{C}$ are the group schemes $\ZZ/p\ZZ$ and~$\mu_p$.

\proof\ The proof is similar to the one of Prop.~5.1, except that we also multiply the group
scheme $G'$ that occurs there with a twisted constant 
group scheme of the type that occurs in part 2 of section~2. More precisely, we multiply $G'$ by
$V=(\ZZ/p\ZZ)^{l-1}$ twisted by the permutation representation
$\varrho:{\rm Gal}(\QQ(\zeta_l^{})/\QQ)\longrightarrow{\rm GL}(V)$. 
The points of~$V(\varrho)$ generate the fidel $\QQ(\zeta^{}_l)$.
Therefore the field $L$ generated by the points of $G'$ satisfies the first condition.

The group scheme $G'$ is killed by~$p$. It is not necessarily an object
of~$\underline{D}$, but it is still an object of~$\underline{C}$.  The second condition follows from the
fact that $G'$ is a $p$-group scheme over~$\ZZ[\hbox{$1\over l$}]$. The third condition is satisfied
because the inertia groups
$I_{\goth l}$ of the primes $\goth l$ over~$l$ act tamely on the points of $G'$, so that the $l$-adic
contribution to the root discriminant of $L$ is of the form $l^{(e-1)/e}$ where $e$ is the ramification
index~$e$ of any of the primes over~$l$. The fourth condition is verified by the theorem of
Abra\v skin-Fontaine. Since the first four conditions are satisfied, the degree
$[L:\QQ(\zeta_p)]$ must be a power of~$p$.
Arguments similar to those used to prove the previous proposition show
that any simple  object of the category~$\underline{C}$ is isomorphic to
$\ZZ/p\ZZ$ or~$\mu_p$ possibly twisted by a character~$\chi$ that is unramified outside~$l$. 
The order of such a character necessarily divides $l-1$ as well as~$p-1$.
 Since
$\zeta_l^{}\in L$, the degree 
$l-1=[\QQ(\zeta_p,\zeta_l):\QQ(\zeta_p)]$ is a power of~$p$. It follows that either $l=2$ or that $p=2$
and 
$l$ is a Fermat prime. In either case we have that  ${\rm gcd}(l-1,p-1)=1$, so that there are
no non-trivial twists of $\ZZ/p\ZZ$ or $\mu_p$ by a character that is unramified
outside~$l$. This proves the proposition.

\beginsection 6. Odlyzko bounds and class field theory.

In this section we prove Theorems~1.1 and Theorem~1.3.  We first deal with Theorem~1.3 and
check the conditions of Prop.3.2 for the pairs of primes
$(l,p)=(2,3)$, $(3,2)$ and~$(5,2)$. Since in each case we have that
${{l^2-1}\over{24}}\not\equiv 0\zmod p$, Proposition~4.1 implies that any extension of
$\mu_p$ by $\ZZ/p\ZZ$ over $\ZZ[{1\over l}]$ is necessarily  trivial.
This leaves us with the second part of the criterion. We check it, case by case, using Proposition~5.2

We recall that the {\it root discriminant} of a number field $K$ of degree $n$ is the $n$-th root of
the absolute value of its absolute discriminant. The second, third and fourth property
of the number field $L$ that
occurs Proposition~5.2 imply that the root discriminant $\delta_L$  satisfies
$$
\delta_L<lp^{1+{1\over{p-1}}}.
$$
We proceed case by case. The strategy is to deduce from the first four properties of the
field $L$ of  Prop.~5.2 that the degree $[L:\QQ(\zeta_p^{})]$ is
necessarily a power of~$p$.

\medskip\noindent{\bf Case $l=2$, $p=3$.\ } The root discriminant $\delta_L$ of the
field~$L$ of Proposition~5.2 satisfies $\delta_L<2 \cdot 3^{3/2}= 10.49\ldots$. Odlyzko's
bounds~[\MA,~p.187] or~[\ODL] imply that
$[L:\QQ]<24$. Therefore the degree of $L$ over  $K=\QQ(\zeta_3,{\root 3\of 2})$  is at most~3. If it
were~2, then the quadratic extension $K\subset L$ would only be ramified at~3. But $K$ does not admit a
quadratic extension that is unramified outside the unique prime over~3. This follows from class field
theory because  the class number of $K$ is~1 and  the multiplicative group $\FF_3^*$ of
the residue field of the unique prime over~3 is generated by the global unit~$-1$. This proves
that $[L:\QQ(\zeta_3)]$ is a power of~3 as required.

\medskip\noindent{\bf Case $l=3$, $p=2$.\ } The root discriminant $\delta_L$ of the
field~$L$ of  Proposition~5.2 satisfies
$\delta_L<3\cdot 2^2= 12$. Odlyzko's bounds imply that $[L:\QQ]<32$.
Therefore the degree of $L$ over  $\QQ(i,\sqrt{-3})=\QQ(\zeta_{12})$  is at most~7. Let
$\pi={\rm Gal}(L/\QQ)$. Since the 2-adic valuation of the root discriminant of
$\QQ(\zeta_{24})$ is not strictly smaller than~2, the field $\QQ(\zeta^{}_{12})$ is the
largest abelian extension of $\QQ$ inside~$L$. It is the fixed field of the commutator
subgroup
$\pi'$ of~$\pi$.  The field $\QQ(\zeta_{12})$ admits no everywhere unramified extension.
Moreover,
$\zeta_3$ generates the multiplicative group $\FF_4^*$ of the residue field of the unique
prime over~2 in $k$ while  the multiplicative group $\FF_9^*$ of the residue field of the
unique prime over~3 is a 2-group. It follows from class field theory  that
$\pi'/\pi''$ is a 2-group. Finally we claim that $\pi''$ is trivial. Since $\#\pi'\le7$,
this is clear when $[\pi':\pi'']=1$ or~4. When $[\pi':\pi'']=2$, the group $\pi''$ has
order at most~3 and is cyclic. Therefore ${\rm Aut}(\pi'')$ is abelian and $\pi'$ is in
the kernel of the homomorphism
$\pi\longrightarrow{\rm Aut}(\pi'')$  induced by conjugation. This implies that $\pi''$ is
contained in the center of~$\pi'$. Since $\pi'/\pi''$ is cyclic, this implies that $\pi'$
is abelian and hence that $\pi''$ is trivial, as required.

\medskip\noindent{\bf Case $l=5$, $p=2$.\ } Let $\pi={\rm Gal}(L/\QQ)$. The root
discriminant
$\delta_L$ of the field~$L$ of Prop.~5.2 satisfies
$\delta_L<5\cdot 2^2= 20$. Odlyzko's bounds imply that $[L:\QQ]<480$ so that the
degree of~$L$ over $\QQ(i,\zeta_{5})=\QQ(\zeta_{20})$ is less than~60. This implies
that $\pi={\rm Gal}(L/\QQ)$ is a solvable group. Since the 2-adic valuation of the root discriminant of
$\QQ(\zeta^{}_{40})$ is not strictly smaller than~2, the field
$\QQ(\zeta_{20})$ is the largest abelian extension of~$\QQ$ inside~$L$. 
Next we study the maximal abelian extension $K$ of $\QQ(\zeta^{}_{20})$ inside~$L$.
It is the fixed field of~$\pi''$.
\smallskip
\noindent {\sl Step 1.} The  field $K$ is contained in the ray class field
$F_{2(1-\zeta_5)}$  of
$\QQ(\zeta_{20}^{})$ of conductor~$\hbox{$2(1-\zeta_5^{})$}$. 

\noindent{\sl Proof.} The class number of
$\QQ(\zeta_{20})$ is~1, the multiplicative group $\FF_{16}^*$ of the residue field of the unique
prime over~2 is generated by  the global unit~$1-\zeta_{20}$ and the multiplicative groups of the residue
fields of the two primes over~5 both have order~4.  By class field theory the field $\QQ(\zeta_{20})$ 
admits therefore no non-trivial unramified odd degree extensions inside~$L$. In particular
${\rm Gal}(K/\QQ(\zeta^{}_{20}))\cong\pi'/\pi''$ has order a power of~2.
It follows that the characters of ${\rm
Gal}(K/\QQ(\zeta_{20}^{}))$ have 2-power order. Their conductors divide $(1-\zeta^{}_5)(1-i)^a$ for some
$a\ge0$. Since the 2-adic valuation of $\delta_L$ is strictly smaller than~2, we have that $a\le 3$.
A character $\chi$ of conductor divisible by $(1-i)^3$ cannot be quadratic. In addition, $\chi^2$ has
conductor divisible by $(1-i)^2$. It follows from the conductor discriminant formula that the 2-adic
valuation of the root discriminant of the field cut out by~$\chi$ is at least 2 which is impossible.
Therefore $a\le 2$ and $K$ is contained in $F_{2(1-\zeta_5^{})}$ as required.
\smallskip
\noindent {\sl Step 2.} The ray class field $F_{2(1-\zeta_5^{})}$ is a biquadratic extension
of~$\QQ(\zeta_{20})$. Its root discriminant is equal to~$5^{7/8}2^{7/4}$.

\noindent{\sl Proof.} The unit group of $\ZZ[\zeta^{}_{20}]$ is generated by
$\zeta_{20}^{}$ and by $1-\zeta_{20}^{a}$  for $a\in(\ZZ/20\ZZ)^*$. Since these global units generate the
multiplicative group of the ring $\ZZ[\zeta_{20}^{}]/(1-\zeta^{}_5)$, class field theory implies that
the ray class field of $\QQ(\zeta_{20}^{})$ of conductor  $1-\zeta_5$ is trivial.
In a similar way,  the global units generate a subgroup of $(\ZZ[\zeta_{20}^{}]/(2))^*$ of index~2.
This implies that the ray class field $F_2$ of $\QQ(\zeta_{20}^{})$ of conductor 
$2$ has degree~2. It is not difficult to see that  $F_2=\QQ(\zeta^{}_{20},\sqrt{\eta})$ where
$\eta=(1+\sqrt{5})/2$. A similar computation shows that the ray class field $F_{2(1-\zeta_5)}$ of 
conductor~$2(1-\zeta^{}_5)$ is a biquadratic extension of~$\QQ(\zeta_{20}^{})$. 
Since the three quadratic characters have conductors $2$, $2(1-\zeta^{}_5)$ and $2(1-\zeta^{}_5)$
respectively, the root discriminant of
$F_{2(1-\zeta_5)}$ is equal to
$5^{7/8}2^{7/4}=13.75\ldots$ as required.
\smallskip
\noindent {\sl Step 3.} The group $\pi''/\pi'''$ is a 2-group.

\noindent{\sl Proof.} The Odlyzko bounds imply that the absolute degree of the Hilbert
class field of
$F_{2(1-\zeta_5)}$ is at most 46. Since $F_{2(1-\zeta_5)}$ has degree 32 over~$\QQ$, it therefore admits
no everywhere unramified extension inside~$L$. The same is true for the three quadratic extensions
of~$\QQ(\zeta_{20})$  contained in $F_{2(1-\zeta_5^{})}$. 
The prime over~2 is totally ramified in 
$F_{2(1-\zeta_5)}$.  Since $F_2= \QQ(\zeta^{}_{20},\sqrt{\eta})$ where
$\eta=(1+\sqrt{5})/2$, the prime over~$5$ is inert in the subfield  $F_2$.
It follows that the residue fields of the primes over~2 and~5 in $F_2$ are isomorphic to
$\FF_{16}$ and~$\FF_{25}$ respectively. The same is true for $F_{2(1-\zeta_5)}$.
The units $1-\zeta^{}_{20}$ and $i\pm\sqrt{\eta}$ of the field $F_2$ 
together with their conjugates generate
a subgroup of the  multiplicative group $\FF_{16}^*\times \FF_{25}^*\times\FF_{25}^*$
of index a power of~2.
Therefore by class field theory neither $F_{2(1-\zeta_5)}$ nor any of the three quadratic
extensions of $\QQ(\zeta_{20})$ contained in it,
admit an odd degree extension inside~$L$. We conclude that 
$\pi''/\pi'''$ is a 2-group as required.
\smallskip

We conclude the proof by showing that $\pi'''$ is trivial. This is clear when
$[\pi':\pi'']\le 2$, because then the order of $\pi''/\pi'''$ is odd by the Burnside basis
theorem. But then $\pi''/\pi'''$ and hence $\pi''$ is trivial. When $[\pi':\pi'']=4$, the
order of~$\pi''$ is at most~14. By Taussky's theorem~[\TAU], the group $\pi''/\pi'''$ is a
{\it cyclic} 2-group. If it is trivial we are done. If not, then $\#\pi'''$ is odd and at
most~7. It follows that $\pi'''$ is cyclic and hence has an abelian automorphism group.
Therefore $\pi''$ is in the kernel of the homomorphism $\pi'\longrightarrow{\rm Aut}(\pi''')$
induced by conjugation. It follows that $\pi'''$ is contained in the center
of~$\pi''$, so that $\pi''$ modulo its center is cyclic. This implies that $\pi''$ is abelian
and hence that $\pi'''$ is trivial, as required.

\medskip
Now we prove Theorem~1.1. Since we already proved Theorem~1.3, there is nothing left to be done 
for the primes~$l=2$, 3 and~5. We deal with the primes~7 and~13 by 
checking the conditions of Proposition~5.1 for the pairs $(l,p)=(7,3)$ and $(13,2)$.
Since in each case we have that ${{l^2-1}\over{24}}\not\equiv 0\zmod p$,
Proposition~4.1 applies and we only need to show that the
simple objects in the category $\underline{D}$ are $\ZZ/p\ZZ$ and $\mu_p$.
This is proved by applying Proposition~5.1.
The first three properties of the number field $L$ that occurs in this proposition imply that
the ramification index of any of the primes over~$l$ in $L$
is at most~$p$ and that the root discriminant $\delta_L$  satisfies
$$
\delta_L<l^{1-{1\over p}}p^{1+{1\over{p-1}}}.
$$
We show  that they imply that
the field $L$ is necessarily of $p$-power degree over~$\QQ(\zeta_p^{})$.

\medskip\noindent{\bf Case $l=7$, $p=3$.\ } In this case the field $F$ of Prop.~5.1 is equal to
$\QQ(\zeta^{}_7+\zeta_7^{-1})$. The root discriminant
$\delta_L$ of the field~$L$ satisfies $\delta_L<3^{3/2}\cdot 7^{2/3}= 19.01\ldots$.
Odlyzko's bounds imply that
$[L:\QQ]<270$ so that the degree of~$L$ over
$K=\QQ(\zeta_3^{},\zeta_7^{}+\zeta_7^{-1},{\root 3 \of 7})$ is at most~14. In particular,
$\pi={\rm Gal}(L/K)$ is a solvable group. Since $7\equiv 1\zmod 3$, the field
$\QQ(\zeta_3^{},{\root 3 \of 7})$ is a cubic extension of conductor $3\cdot 7$ of
$\QQ(\zeta_3^{})$. It follows from the conductor discriminant formula and the relative
discriminant formula that the root discriminant of $K$ is equal to
$3^{7/6}7^{2/3}=13.18\ldots$. By Odlyzko's bounds, any unramified extension of~$K$ has relative
degree at most~$42/18$. This implies that $K$ admits at most an unramified quadratic extension~$H$. By
group theory such a field $H$ would be  the composite of~$K$ with an unramified quadratic extension
of~$\QQ(\zeta_3)$. This implies that $K=H$. There lies only one prime $\pp$ over~3 in~$K$. Its residue
field has 27 elements. It is easily seen that its unit group is generated by the global units~$-1$ and
$\zeta_7^{}+\zeta_7^{-1}$. It follows therefore from class field theory that there is no abelian
extension of $K$ inside
$L$ that is at most tamely ramified at the prime over~$3$. This implies that $\pi/\pi'$ is a 3-group. 

If $[\pi:\pi']=1$ or~9, we are done. Suppose that $[\pi:\pi']=3$  and let $\pp^a$ be the conductor
of the corresponding cubic extension $K'$ of~$K$. By the conductor discrimiannt formula, the root
discriminant of $K'$ is equal to $3^{a/9+7/6}7^{2/3}$. Since it is strictly smaller than 
$3^{3/2}7^{2/3}$, we must have that $a=2$. This implies that the root dsicriminant of $K'$ is
$3^{25/18}7^{2/3}=16.82\ldots$. Odlyzko's bounds imply that any unramified extension of $K'$ has
relative degree smaller than~$120/36<4$. The field $K'$ does not admit a quadratic extension that is at
most tamely ramified at the prime over~3 for the same reasons as above. It follows that
$[L:K']$ and hence $[L:\QQ(\zeta_3^{})]$ is a power of~3, as required.

\medskip\noindent{\bf Case $l=13$, $p=2$.\ }  The root discriminant $\delta_L$ of the
extension~$L$ satisfies $\delta_L<2^2\cdot 13^{1/2}= 14.422\ldots$. Odlyzko's bounds imply
that $[L:\QQ]<60$ so that the degree of~$L$ over $\QQ(i,\sqrt{13})$ is at most~14. Let
$\pi={\rm Gal}(L/\QQ)$. Since the root discriminant of $\QQ(\zeta_8,\sqrt{13})$ is equal
to~$4\sqrt{13}>\delta_L$, the largest abelian extension of $\QQ$ inside~$L$ is $\QQ(i,\sqrt{13})$. This
implies that $\pi/\pi'$ is isomorphic to the Klein four group. The class number of
$\QQ(i,\sqrt{13})$ is~1 and since $\FF_4^*$ is generated by the global
unit $\eta=(3+\sqrt{13})/2$, the ray class field of $\QQ(i,\sqrt{13})$ of conductor~$(1+i)$ is trivial.
By class field theory
$\pi'/\pi''$ is a 2-group. Suppose that $\chi$ is a quadratic character of
$\pi'/\pi''$. Then ${\rm cond}(\chi)=(1+i)^a$ for some $a\le 2$. Since $v_2(\delta_L)<2$,
we have $a\le 3$. Since there are no quadratic characters of conductor $(1+i)^3$ we
must have that $a=2$. Since the unit $i$ is not congruent to~1 modulo $(1+i)^2$, the ray class field
of conductor $(1+i)^2$ has degree~2 over $\QQ(i,\sqrt{13})$. It is the field generated by
$\sqrt{\eta}$. Any larger abelian extension of $\QQ(i,\sqrt{13})$ inside $L$ has relative
discriminant at least $(i+1)^{2+3+3}$ and hence root discriminant at least $4\sqrt{13}$. This is
impossible and hence $\pi'/\pi''$ has order at most~2. 

If $\#\pi'/\pi''=1$, it follows that $\pi'=1$ and hence that  $\pi$ is a 2-group. If
$\#\pi'/\pi''=2$, we know that $K=\QQ(i,\sqrt{13},\sqrt{\eta})$ is the fixed field of
$\pi''$ and that $\pi''/\pi'''$ has odd order. However, $K$ admits no odd degree abelian extensions
inside~$L$. Indeed, $\delta_K$ is equal to
$2^{3/2}\sqrt{13}=10.198\ldots$ so that Odlyzko's bounds imply that the absolute degree of its Hilbert
class field is at most 22. Since
$[K:\QQ]=8$, this implies that the class number of $K$ is at most~2. Since the residue field of the
unique prime
$\goth p$ over~2 in $F$ is equal to the residue field~$\FF_4$ of the prime $i+1$ in
$\QQ(i,\sqrt{13})$, the ray class field of
$K$ of conductor~$\goth p$ is equal to $K$ itself. This implies that $K=L$ and that
$[L:\QQ]$ is a power of~2 as required.

\beginsection 7. The case $l=11$.

In this section we prove Theorem~1.2. The modular curve $X_0(11)$ has genus~1 and is given
by the equation
$$
Y^2 +Y = X^3 -X^2 -10X -20.
$$ 
Its Jacobian $E=J_0(11)$ is semi-stable and has good reduction at every prime except at~$l=11$. We take
$p=2$ and study $p$-group schemes over~$\ZZ[\hbox{$1\over{11}$}]$ in the
category~$\underline{D}$ introduced in section~2. The 2-torsion points
$E[2]$ of $E$ form a 2-group scheme  that is an object in~$\underline{D}$ .
The natural map
${\rm Gal}(\overline{\QQ}/\QQ)\longrightarrow{\rm Aut}(E[2](\overline{\QQ}))$ is surjective.
Indeed,
the points in $E[2]$ generate the sextic field $F=\QQ(\sqrt{-11},\alpha)$ where $\alpha$
satisfies the equation $\alpha^3+\alpha^2+\alpha-1=0$. This implies that $E[2]$ is a {\it
simple} object of~$\underline{D}$. It is isomorphic to its Cartier dual. Since the elliptic
curve $E$ is supersingular modulo~2, the group scheme~$E[2]$ is local-local (i.e. local with
local Cartier dual) at~2.

\proclaim Proposition 7.1. The simple objects in the category~$\underline{D}$ are 
the group schemes $\ZZ/2\ZZ$, $\mu_2$ and~$E[2]$.

\proof\ We modify the proof of Proposition~5.1. Let $G$ be simple and let $G'$ be the
product of
$G$ by $E[2]$ and by the Katz-Mazur group schemes $G_{\varepsilon}$ for the units $\varepsilon=-1$
and~$11$ of~$\ZZ[{1\over{11}}]$. Then $G'$ is killed by~2 and is again an object of~$\underline{D}$.
Therefore the root discriminant $\delta_L$ of the extension $L$ generated by the points of~$G'$
satisfies
$\delta_L<4\sqrt{11}=13.266\ldots$. Odlyzko's bounds imply that $[L:\QQ]<44$. We have the
inclusions
$$
\QQ\quad\mathop{\subset}\limits_4\quad\QQ(i,\sqrt{-11})\quad\mathop{\subset}\limits_3 \quad F(i)
\quad\mathop{\subset}\limits_{\le3}\quad L.
$$
It follows that $[L:F(i)]\le3$.  If $[L:F(i)]=3$, the field  $L$ is abelian
over~$\QQ(i,\sqrt{-11})$. Since the class number of $\QQ(i,\sqrt{-11})$ is~1 and since the ray
class field of conductor~2 of $\QQ(i,\sqrt{-11})$ is $F(i)$, class field theory implies that $[L:F(i)]$
cannot be equal to~3 and must therefore be 1~or~2. In either case, ${\rm Gal}(L/F(i))$ and hence ${\rm
Gal}(L/F)$ is a 2-group and hence fixes the points of the simple 2-group scheme~$G$. Therefore
${\rm Gal}(\overline{\QQ}/\QQ)$ acts on
$G(\overline{\QQ})$ through ${\rm Gal}(F/\QQ)\cong S_3$. It follows from the structure of
simple
$\FF_2[S_3]$-modules that the Galois module
$G(\overline{\QQ})$ is either a group of order~2 with trivial action or is isomorphic
to~$E[2](\overline{\QQ})$.

If $G$ has order~2, it follows that $G\cong\ZZ/2\ZZ$ or~$\mu_2$. If 
$G(\overline{\QQ})\cong E[2](\overline{\QQ})$, the inertia group $I_2$ of the prime over~2
acts irreducibly, so $G$ is local-local. By
Raynaud's theorem~[\RAY,~section 3.3.5] the group scheme $G$ is therefore determined by its
Galois module and hence we have that $G\cong E[2]$ as group schemes over~$\ZZ_2$. This leads to an
isomorphism of local Galois modules $G(\overline{\QQ}_2)\cong E(\overline{\QQ}_2)$
Since the isomorphisms $G(\overline{\QQ})\cong E[2](\overline{\QQ})$ and
$G(\overline{\QQ}_2)\cong E[2](\overline{\QQ}_2)$ are unique, they are compatible.
It follows from the equivalence of categories of~[\SCH,~Prop.2.4] that
$G\cong E[2]$ over the ring~$\ZZ[{1\over{11}}]$. This proves the Proposition.

\medskip
Next we study extensions in~$\underline{D}$ of the simple group schemes with one another.
Note that the group scheme $E[4]$ fits in the non-split exact sequence 
$$
0\longrightarrow\quad E[2]\quad\longrightarrow \quad E[4]\quad\longrightarrow\quad
E[2]\quad\longrightarrow 0
$$
and is an object of~$\underline{D}$, because $E$ has semi-stable reduction at~11.

\proclaim Proposition 7.2. Over the ring~$\ZZ[{1\over{11}}]$ we have the following.
\item{(i)} the groups ${\rm Ext}^1_{\ZZ[{1\over{11}}]}(\mu_2,\ZZ/2\ZZ)$, ${\rm
Ext}_{\ZZ[{1\over{11}}]}^1(E[2],\ZZ/2\ZZ)$ and
${\rm Ext}^1_{\ZZ[{1\over{11}}]}(\mu_2,E[2])$ are all trivial.
\item{(ii)} The group of extensions ${\rm Ext}^1_{\underline{D}}(E[2],E[2])$ has dimension~1
over~$\FF_2$. It is generated by the extension~$E[4]$ above.

\proof\ {\sl (i)} Any extension
$$
0\longrightarrow\ZZ/2\ZZ\longrightarrow G\longrightarrow\mu_2\longrightarrow 0
$$
is split over $\ZZ_2$ by the connected component. Therefore $G$ is killed by 2 and its 2-adic Galois
representation is trivial. It follows that ${\rm Gal}(\overline{\QQ}/\QQ)$ acts on $G(\overline{\QQ})$
through a character $\chi$ that is at most ramified at~$11\cdot\infty$. There is exactly one
such
$\chi$ that is non-trivial. The prime 2 is {\it inert} in the corresponding
field~$\QQ(\sqrt{-11})$. Therefore
$\chi=1$ and $G$ is generically split. The Mayer-Vietoris sequence~[\SCH,~Cor.2.4] implies then that $G$
is split over~$\ZZ[{1\over{11}}]$, as required.

The proof that every extension  $0\longrightarrow\ZZ/2\ZZ\longrightarrow
G\longrightarrow E[2]\longrightarrow 0$ splits is entirely similar and the fact that ${\rm
Ext}^1_{\ZZ[{1\over{11}}]}(\mu_2,E[2])$ vanishes follows  by Cartier duality.

\smallskip\noindent
{\sl (ii)\ } Since the Galois representation ${\rm Gal}(\overline{\QQ}/\QQ)\longrightarrow{\rm
Aut}(E[2](\overline{\QQ}))$ is surjective, the only  Galois equivariant endomorphisms of the
group scheme $E[2]$ are scalar multiplications. This is true over any of the rings $\ZZ_2$,
$\QQ_2$, $\ZZ[\hbox{$1\over{11}$}]$ and~$\ZZ[\hbox{$1\over{22}$}]$.

It follows from the Mayer-Vietoris sequence
that there is an exact sequence
$$
0\!\rightarrow{\rm Ext}^1_{\ZZ[{1\over{11}}]}(E[2],E[2])\longrightarrow {\rm
Ext}^1_{\ZZ_2}(E[2],E[2])\times {\rm Ext}^1_{\ZZ[{1\over{22}}]}(E[2],E[2])
\longrightarrow {\rm Ext}^1_{\QQ_2}(E[2],E[2])\rlap.
$$
The proof proceeds in two steps.
\medskip
\noindent {\it Step 1.} {\sl  The group ${\rm Ext}^1_{\underline{D}}(E[2],E[2])$ is generated by $E[4]$
and by the extensions in~$\underline{D}$ that are killed by~2.}

\smallskip\noindent Let $\Gamma={\rm Gal}(\overline{\QQ}/\QQ)$ and let ${\rm Ext}^q_{\rm
ab}$ denote the $q$-th Ext-group in the category of abelian groups.  From the spectral sequence 
$$H^p(\Gamma,{\rm Ext}^q_{\rm ab}(E[2],E[2]))\quad\Longrightarrow\quad
{\rm Ext}^{p+q}_{\QQ}(E[2],E[2])
$$ 
we deduce the following  commutative diagram with exact
rows
$$
\def\normalbaselines{\baselineskip20pt\lineskip3pt 
\lineskiplimit3pt}
\matrix{
0&\longrightarrow&{\rm Ext}^1_{\underline{D}}(E[2],E[2])_2&\longrightarrow&
{\rm Ext}^1_{\underline{D}}(E[2],E[2])&\longrightarrow&{\rm cok}&\!\!\!\!\longrightarrow&0\cr
&&\Big\downarrow&&\Big\downarrow&&\Big\downarrow&&\cr
0&\longrightarrow&{\rm Ext}^1_{\QQ}(E[2],E[2])_2&\longrightarrow&
{\rm Ext}^1_{\QQ}(E[2],E[2])&\longrightarrow&{\rm Ext}_{\rm
ab}^1(E[2],E[2])^{\Gamma}&&\cr}
$$
Here the index `$2$' means `annihilated by~2'.
Any extension $G$ of $E[2]$ by $E[2]$ over $\ZZ[\hbox{$1\over{11}$}]$ that is killed by~2
over~$\QQ$, is itself killed by~2. Therefore the left hand square is Cartesian. It
follows that the rightmost vertical arrow is injective. Since the
$\Gamma$-modules
${\rm Ext}_{\rm ab}^1(E[2],E[2])$ and ${\rm Hom}_{\rm ab}(E[2],E[2])$ 
are dual to one another, the orders of
${\rm Ext}_{\rm ab}^1(E[2],E[2])^{\Gamma}$ and ${\rm Hom}_{\rm ab}(E[2],E[2])^{\Gamma}
={\rm Hom}^{}_{\Gamma}(E[2],E[2])$ are equal. 
Since the  Galois represention on the points of~$E[2]$
is surjective, the latter group consists of the scalar matrices and has order~2. It
follows that the index of ${\rm Ext}^1_{\underline{D}}(E[2],E[2])_2$ in ${\rm
Ext}^1_{\underline{D}}(E[2],E[2])$ is at most~2. It is exactly 2, because of the existence of
the group scheme~$E[4]$.

\medskip
\noindent {\it Step 2.} {\sl Any extension {\it in the category $\underline{D}$} of $E[2]$ by~$E[2]$
over~$\ZZ[{1\over{11}}]$ that is annihilated by~2, is trivial.} 

\smallskip\noindent Using Odlyzko's bounds one shows that the field $F=\QQ(\sqrt{-11},\alpha)$ has
class number~1. Let $\pi\in F$ be a prime over~2. We have that $(\pi)^3=(2)$. Let $G$ be
an extension of $E[2]$ by~$E[2]$ that is on object in~$\underline{D}$ and that is annihilated
by~2. The field extension
$L$ generated by the points of~$G$ is of exponent~2 over~$F$ and is at most ramified at the primes
over~2 and~11. Since $G$ is an object in $\underline{D}$ of exponent~2, the inertia groups $I_{11}
\subset {\rm Gal}(L/\QQ)$ of the primes over~11, have order at most~2. Since the ramification index is
equal to~2 in the extension $F$ of~$\QQ$, this implies that
$L$ is actually {\it unramified} over~$F$ at the primes over~11.
By~[\SCH,~Prop.6.4],  the field $L$ is a biquadratic extension  of~$F$ of  conductor
dividing~$\pi^2$. Since the global unit $\alpha$ generates the group $(1+(\pi))/(1+(\pi^2))$ and since $F$
admits no non-trivial unramified extensions, class field theory implies that $F=L$. Therefore
the extension
$G$ is split as an extension of Galois modules. It follows from~[\SCH,~Prop.6.4] (or rather the
proof of its part~{\sl (ii)}) that then the extension is necessarily locally trivial. The
Mayer-Vietoris sequence  implies then that $G$ is split, as required.

This proves~{\sl (ii)}.
\medskip
Note that the space ${\rm Ext}^1_{\ZZ[{1\over{11}}]}(E[2],E[2])$ of {\it all} extensions
of~$E[2]$ by $E[2]$ does {\it not} have dimension~1. Indeed, consider the 
elliptic curve $E'$ given by the Weierstrass equation $Y^2+Y=X^3-X^2-7X+10$. It is the curve 121D
in~[\ANT,~p.97] of conductor~$11^2$. Over the ring $\ZZ[{1\over{11}}]$ the subgroup scheme
$E'[2]$ of 2-torsion points is isomorphic to the subgroup scheme of 2-torsion points of the
semi-stable abelian variety~$E=J_0(11)$. This follows from the fact that the points of each
of the two group schemes generate the field  $F=\QQ(\sqrt{-11},\alpha)$  with
$\alpha^3+\alpha^2+\alpha-1=0$. A theorem of Raynaud's~[\RAY,~section 3.3.5] implies then
that $E'[2]\cong E[2]$ over~$\ZZ_2$ and it follows~[\SCH,~Prop.2.4] that the same is
true over~$\ZZ[{1\over{11}}]$. Therefore $E'[2]$ is an object in~$\underline{D}$.  

This is in agreement with the fact that in the short table below~[\ANT,~Table~3], the
coefficients $a(p)$ and $b(p)$ of the $L$-functions of
$E'$ and $E$ are congruent modulo~2. The coefficients are visibly {\it not} congruent
modulo~4. This implies that the group
scheme $E'[4]$ is {\it not} isomorphic to
$E[4]$. Since ${\rm Ext}^1_{\underline{D}}(E[2],E[2])$ is 1-dimensional, we conclude
that $E'[4]$ is not an object of~$\underline{D}$ and hence that the dimension of 
${\rm Ext}^1(E[2],E[2])$ is at least~2. Being an extension of~$E'[2]$ by~$E'[2]$, the group
scheme $E'[4]$ is still an object of the category~$\underline{C}$.

\vbox{
\medskip
\centerline {\vbox {\offinterlineskip
\hrule\halign{&\vrule#&\strut\quad\hfil#\quad\cr
height3pt
&\omit&&\omit&&\omit&&\omit&&\omit&\cr
&$p$&&2&&3&&5&&7&\cr
height3pt
&\omit&&\omit&&\omit&&\omit&&\omit&\cr
\noalign{\hrule}
height2pt
&\omit&&\omit&&\omit&&\omit&&\omit&\cr
&$a_{p}^{}$&&$-2$&&$-1$&&1&&$-2$&\cr
height2pt
&\omit&&\omit&&\omit&&\omit&&\omit&\cr
&$b_{p}^{}$&&$0$&&$-1$&&$-3$&&$0$&\cr
height3pt
&\omit&&\omit&&\omit&&\omit&&\omit&\cr
}\hrule}}}

\medskip\noindent{\bf Proof of Theorem~1.1.} Let $A$ be a semi-stable abelian variety over~$\QQ$ with
good reduction outside~11. We show that the subgroup scheme $A[2]$ of 2-torsion points admits a
filtration with successive subquotients isomorphic to~$E[2]$. Consider an arbitrary filtration
of~$A[2]$ with simple subquotients. By Prop.~7.2~{\sl (i)} we obtain a filtration of the form
$$
0\quad\subset\quad G_1\quad\subset\quad G_2\quad\subset\quad A[2]
$$
with $G_1$ an extension of group schemes isomorphic to $\mu_2$, with $A[2]/G_2$ an extension of
copies of $\ZZ/2\ZZ$ and
$G_2/G_1$ admitting a filtration with copies of~$E[2]$.
We claim that actually, there are no subquotiens isomorphic to $\ZZ/2\ZZ$ or $\mu_2$ in this
filtration. Indeed, suppose that there is a subquotient isomorphic to $\ZZ/2\ZZ$. Using
Prop.~7.2~{\sl (i)} to `move the subquotients that are isomorphic to $\ZZ/2\ZZ$ to the right', we find
that for each
$n\ge 1$  there is an exact sequence
$$ 
0\longrightarrow H_n \longrightarrow A[2^n]\longrightarrow C_n\longrightarrow 0
$$ 
with $C_n$ admitting a filtration of length~$n$ with subquotients isomorphic to~$\ZZ/2\ZZ$. 
The fundamental group of $\ZZ[{1\over{11}}]$ acts on $C_n(\overline{\QQ})$ through a
2-group~$P$. Since the maximal abelian 2-extension of~$\QQ$ that is unramified outside~11 is
the quadratic field~$\QQ(\sqrt{-11})$, the groups $P/P'$ and hence~$P$ are cyclic and the
group schemes 
$C_n$ become constant over the ring $\ZZ[{{1+\sqrt{-11}}\over 2},{1\over{11}}]$. Now choose a
non-zero prime
$\pp$ of the ring $\ZZ[{{1+\sqrt{-11}}\over 2},{1\over{11}}]$. The  abelian varieties $A/H_n$ are all
isogenous to~$A$ and have therefore the same number of points modulo~$\pp$.
On the other hand, they have at least~$2^n$ rational points. This leads to a contradiction when
$n\rightarrow\infty$. Therefore there are no subquotients isomorphic to $\ZZ/2\ZZ$ in the
filtration. By Cartier duality there are none isomorphic to 
$\mu_2$ either. 

It follows that $A[2^n]$ is filtered with group schemes isomorphic to~$E[2]$. Since $E$ has no
complex multiplication,  Tate's theorem~[\TDB] implies that  the endomorphism ring of the 2-divisible
group
$G$ of $E$ is isomorphic to~$\ZZ_2$. An application of Theorem~8.3 of the next section to the
2-divisible group
$H$ of~$A$ over the ring $O=\ZZ[{1\over{11}}]$ shows that the 2-divisible groups of
$A$ and $E^g$ are isomorphic. By Faltings' theorem~[\FAL] this implies that $A$ is isogenous to~$E^g$
as required.

\beginsection 8. $p$-divisible groups.

The main result of this section is Theorem~8.3. It is a general result about
$p$-divisible groups, used in section~7 of this paper. See~[\OPDV, \SCHO] for related statements.

Let $O$ be a Noetherian domain of characteristic~0, let $p$ be a prime and let $\underline{D}$ be a full subcategory
of the category of $p$-group schemes over~$O$ that is closed under taking products, closed flat subgroup schemes and
quotients by closed flat subgroup schemes. In section~7, Theorem~8.3 is applied
to the category~$\underline{D}$ that was introduced in section~2.

We first prove a Lemma.

\proclaim Lemma 8.1. Let  $G=\{G_n\}$ be a $p$-divisible group over~$O$. Suppose that $R={\rm End}(G)$
is a discrete valuation ring with uniformizer~$\pi$ and residue field $k=R/\pi R$.
Suppose that
\smallskip
\item{--} every group scheme $G_n$ is an  object in the category~$\underline{D}$;
\smallskip
\item{--} the  map
$$
{\rm
Hom}_O(\gp,\gp)\mathop{\longrightarrow}\limits^{\delta} {\rm
Ext}^1_{\underline{D}}(\gp,\gp)
$$ 
associated to the exact sequence $0\rightarrow \gp\rightarrow G[\pi^2]\rightarrow\gp\rightarrow 0$ is
an isomorphism of 1-dimensional
$k$-vector spaces. 
\smallskip\noindent Then 
\smallskip 
\item{(i)}
For all $m,m'\ge 1$ the canonical map
$$
(R/\pi^{m'}R)[\pi^m]\longrightarrow {\rm Hom}_O(G[\pi^{m'}],G[\pi^m])
$$
is an isomorphism
.
\smallskip 
\item{(ii)} For every $m\ge 1$ the group ${\rm Ext}^1_{\underline{D}}(\gp,G[\pi^m])$ is
generated by the extension
$$
0\longrightarrow\quad G[\pi^m]\quad\longrightarrow\quad G[\pi^{m+1}] \quad\mathop{\longrightarrow}\limits^{\pi^m}
\quad G[\pi] \longrightarrow  0.
$$

\smallskip\noindent{\bf Proof.} Suppose that $0\not=f\in R$ has~$G[\pi^k]$ in its kernel.
Let $\overline{F}$ be an
algebraic closure of the quotient field $F$  of~$O$ and 
let $T$ denote  the injective limit  of the groups of points~$G_n(\overline{F})$. 
For $n\ge 1$ let $f_n$ denote the induced map $G_n(\overline{F})\longrightarrow G_n(\overline{F})$
and let
$f_{\rm ind}:T\rightarrow T$ be the
 limit. We have a commutative diagram with exact
rows
$$
\def\normalbaselines{\baselineskip20pt\lineskip3pt 
\lineskiplimit3pt}
\matrix{0&\longrightarrow&G[\pi^k](\overline{F})&\longrightarrow&T&
\mathop{\longrightarrow}\limits^{\pi^k}&T&\longrightarrow&0\cr
&&\Big\downarrow\rlap{\hbox{ $\scriptstyle \subset$}}&&\Big\Vert&&&&\cr
0&\longrightarrow&\ker(f_{\rm
ind})&\longrightarrow&T&\mathop{\longrightarrow}\limits^{f_{\rm ind}}&T&\longrightarrow&0.}
$$ 
Since
$G[\pi^n](\overline{F})\subset {\rm ker}(f_{\rm ind})$, there exists a Galois equivariant
homomorphism
$h:T\longrightarrow T$ for which $h\cdot\pi^k= f_{\rm ind}$. The map $h$ is the limit of
a system of compatible homomorphisms $h_n:G_n(\overline{F})\longrightarrow G_n(\overline{F})$. Let
$T_pG$ denote the Tate module of~$G$ and let $h_{\rm inv}:T_pG\longrightarrow T_pG$ denote the inverse
limit of the
$h_n$. Writing
$f_{\rm inv}$ for the inverse limit of the~$f_n$, we have that 
$h_{\rm inv}\cdot
\pi^k=f_{\rm inv}$. It follows from Tate's Theorem~[\TDB] that $f$ is divisible
by~$\pi^k$ in~$R$. This shows that the maps in~{\sl (i)} are injective.

To prove surjectivity of the maps in~{\sl (i)}, we observe that both sides are finite groups and we
count their orders. The left hand side has order $q^{{\rm min}(m,m')}$ where $q=\#k$. It follows from the
multiplicativity of orders in exact sequences that
$$\#{\rm Hom}_O(G[\pi^m],G[\pi^{m'}])\le\cases{\#{\rm Hom}_O(\gp,G[\pi^{m'}])^m,\cr
\#{\rm Hom}_O(G[\pi^m],\gp)^{m'}.\cr}
$$
Therefore it suffices to show for all $m\ge 1$ that ${\rm Hom}_O(\gp,G[\pi^m])$ and ${\rm
Hom}_O(G[\pi^m],\gp)$ both are 1-dimensional $k$-vector spaces. By assumption this is so for $m=1$.
It is  not clear for $m>1$, because we do not know a priori that an arbitrary morphism of group
schemes $f:G[\pi]\longrightarrow G[\pi^m]$ commutes with~$\pi$. We show contemporarily that the natural
maps
${\rm Hom}_O(\gp,\gp)\longrightarrow {\rm Hom}_O(\gp,G[\pi^m])$ are isomorphisms and that the natural
maps
${\rm Ext}^1_O(\gp,G[\pi^{m+1}])\longrightarrow  {\rm Ext}^1_O(\gp,G[\pi^m])$ are injections. Since by
assumption ${\rm Ext}^1_O(\gp,\gp)$ has dimension~1, this shows that ${\rm dim}\,{\rm
Ext}^1_O(\gp,G[\pi^{m+1}])=1$ for all~$m\ge 1$. The fact that the  natural homomorphisms 
${\rm Hom}_O(G[\pi^m],\gp)\longrightarrow {\rm Hom}_O(\gp,\gp)$ are isomorphism follows in a similar
way. It can also be deduced from Cartier duality. 

To the
infinite commutative diagram with exact columns
$$
\def\normalbaselines{\baselineskip20pt\lineskip3pt 
\lineskiplimit3pt}
\matrix{0&&0&&0&&\cr
\Big\downarrow&&\Big\downarrow&&\Big\downarrow&&\cr
\gp&=&\gp&=&\gp&=&\cdots\cr
\Big\downarrow&&\Big\downarrow&&\Big\downarrow&\cr
G[\pi^2]&\hookrightarrow&G[\pi^3]&\mathop{\hookrightarrow}\limits&G[\pi^4]&
\hookrightarrow&\cdots\cr
\Big\downarrow\rlap{\hbox{ $\scriptstyle \pi$}}&&\Big\downarrow\rlap{\hbox{ $\scriptstyle
\pi$}}&&\Big\downarrow\rlap{\hbox{ $\scriptstyle \pi$}}&&\cr
\gp&\hookrightarrow&G[\pi^2]&\mathop{\hookrightarrow}\limits&G[\pi^3]&
\hookrightarrow&\cdots\cr
\Big\downarrow&&\Big\downarrow&&\Big\downarrow&&\cr
0&&0&&0&&\cr}
$$
we apply the functor ${\rm Hom}_O(G[\pi],-)$ 
and form the associated long sequences of ${\rm Ext}^1_{\underline{D}}$-groups.
These exist by the remarks in section~2. The resulting diagram has exact columns.

The isomorphisms and zero maps in the exact first column are a direct consequence of the
assumptions. The statement for $m=2$ follows at once. The map $g_2$ is the same map as the first morphism
in the first column. So $g_2$ is  an isomorphism and it follows at once that $f_3$ 
is an isomorphism as well. This implies that 
that  $f_2$ and~$f_4$ are both zero and that $f_1$ is an isomorphism.  Finally
$f_5$ is injective and $g_1$ is an isomorphism. This implies the statement for~$m=3$. Note that the map
$g_3$ is the same map as~$g_1$. Now one proceeds inductively.

This proves the lemma.
$$
\def\normalbaselines{\baselineskip20pt\lineskip3pt 
\lineskiplimit3pt}
\matrix{0&&0&&0&&\cr
\Big\downarrow&&\Big\downarrow&&\Big\downarrow&&\cr
{\rm Hom}_O(\gp,\gp)&=&{\rm Hom}_O(\gp,\gp)&=&{\rm
Hom}_O(\gp,\gp)&=&\cdots\cr 
\Big\downarrow\rlap{$\vcenter{\hbox{$\scriptstyle
\cong$}}$}&&\Big\downarrow\rlap{$\vcenter{\hbox{$\scriptstyle
f_1$}}$}&&\Big\downarrow&\cr
{\rm
Hom}_O(\gp,G[\pi^2])&\mathop{\rightarrow}\limits^{g_1}&{\rm
Hom}_O(\gp,G[\pi^3])&\mathop{\rightarrow}\limits&{\rm
Hom}_O(\gp,G[\pi^4])&
\rightarrow&\cdots\cr
\Big\downarrow\rlap{$\vcenter{\hbox{$\scriptstyle
0$}}$}&&\Big\downarrow\rlap{$\vcenter{\hbox{$\scriptstyle
f_2$}}$}&&\Big\downarrow&&\cr
{\rm
Hom}_O(\gp,\gp)&\mathop{\rightarrow}\limits^{g_2}&{\rm
Hom}_O(\gp,G[\pi^2])&\mathop{\rightarrow}\limits^{g_3}&{\rm
Hom}_O(\gp,G[\pi^3])&
\rightarrow&\cdots\cr
\Big\downarrow\rlap{$\vcenter{\hbox{$\scriptstyle
\cong$}}$}&&\Big\downarrow\rlap{$\vcenter{\hbox{$\scriptstyle
f_3$}}$}&&\Big\downarrow&&\cr
{\rm Ext}^1_{\underline{D}}(\gp,\gp)&=&{\rm Ext}^1_{\underline{D}}(\gp,\gp)&=&{\rm
Ext}^1_{\underline{D}}(\gp,\gp)&=&\cdots\cr
\Big\downarrow\rlap{$\vcenter{\hbox{$\scriptstyle
0$}}$}&&\Big\downarrow\rlap{$\vcenter{\hbox{$\scriptstyle
f_4$}}$}&&\Big\downarrow&\cr
{\rm Ext}^1_{\underline{D}}(\gp,G[\pi^2])&\rightarrow&{\rm
Ext}^1_{\underline{D}}(\gp,G[\pi^3])&\mathop{\rightarrow}\limits&{\rm
Ext}^1_{\underline{D}}(\gp,G[\pi^4])&
\rightarrow&\cdots\cr
\Big\downarrow&&\Big\downarrow\rlap{$\vcenter{\hbox{$\scriptstyle
f_5$}}$}&&\Big\downarrow&&\cr
{\rm Ext}^1_{\underline{D}}(\gp,\gp)&\rightarrow&{\rm
Ext}^1_{\underline{D}}(\gp,G[\pi^2])&\mathop{\rightarrow}\limits&{\rm
Ext}^1_{\underline{D}}(\gp,G[\pi^3])&
\rightarrow&\cdots\cr}
$$

\proclaim Corrolary 8.2.  Let  the ring $O$ and the category $\underline{D}$ be as above and
let $G=\{G_n\}$ be a
$p$-divisible group over~$O$. Suppose that
$R={\rm End}(G)$ is a discrete valuation ring with uniformizer~$\pi$ and residue field $k=R/\pi R$.
Suppose that the conditions  of Lemma~8.1 are satisfied. Then every $p$-group scheme in~$\underline{D}$
that admits a filtration  with closed flat subgroup schemes and successive subquotients isomorphic
to~$\gp$, is isomorphic to a group scheme of the shape
$$
\mathop{\oplus}\limits_{i=1}^rG[\pi^{n_i}].
$$

\smallskip\noindent{\bf Proof.} Let $J$ be a such a group scheme. Proceeding by induction
we may assume that there is an exact sequence
$$
0\longrightarrow \mathop{\oplus}\limits_{i=1}^rG[\pi^{n_i}]\longrightarrow J
\longrightarrow \gp \longrightarrow 0.
$$
The class of this extension is an element in the group
$${\rm
Ext}^1_{\underline{D}}(\gp,\mathop{\oplus}\limits_{i=1}^rG[\pi^{n_i}])\quad\cong\quad
\mathop\oplus\limits_{i=1}^r {\rm
Ext}^1_{\underline{D}}(\gp,G[\pi^{n_i}])\quad\cong\quad \FF_2^r
$$ 
The second isomorphism follows from Lemma~8.1~{\sl (ii)}.
The extensions of the form
$$
0\longrightarrow G[\pi^{n_j}]\times W\longrightarrow 
G[\pi^{n_j+1}]\times W \longrightarrow \gp \longrightarrow 0.
$$
with $W=\mathop{\oplus}\limits_{i\not=j}G[\pi^{n_i}]$ form a basis for the vector space $\FF_2^r$.
These extensions have the required shape. Therefore, if we show that the 
Baer sum of two extensions of the right shape is
again of the right shape, we are done. To do this, it suffices to show that kernels and cokernels of
morphism between group schemes of the shape 
$\mathop{\oplus}\limits_{i=1}^rG[\pi^{n_i}]$ are again of that shape. 
By duality it is enough to show this for kernels.

Let therefore
$$
\mathop{\oplus}\limits_{i=1}^rG[\pi^{n_i}]\quad\mathop{\longrightarrow}\limits^g 
\quad\mathop{\oplus}\limits_{j=1}^sG[\pi^{m_j}]
$$
be a homomorphism of group schemes or, equivalently, of sheaves for the fppf topology on~$O$. By
Lemma~8.1~{\sl (i)} there are endomorphisms $f_{ij}\in R={\rm End}_O(G)$ that induce~$g$.
Hence the kernel $K$ of~$g$ is isomorphic to the kernel of the restriction of
the homomorphism
$$
\pmatrix{f_{11}&\cdots&f_{r1}\cr
\vdots&&\vdots\cr  f_{1s}&\cdots&f_{rs}\cr}:G^r \quad\longrightarrow\quad G^s.
$$
to the subgroup scheme $\mathop{\oplus}\limits_{i=1}^rG[\pi^{n_i}]$ of~$G$.
Consider the commutative diagram 
$$
\def\normalbaselines{\baselineskip20pt\lineskip3pt 
\lineskiplimit3pt}
\matrix{&&0&&0&&&&\cr
&&\downarrow&&\downarrow&&&&\cr
&&K&&K_1&&&&\cr
&&\Big\downarrow&&\Big\downarrow&&&&\cr
0&\longrightarrow&\mathop\oplus\limits_{i=1}^rG[\pi^{n_i}]&
\longrightarrow&G^r&
\mathop{\longrightarrow}\limits^{\hbox{$ \scriptstyle {\pmatrix{\pi^{n_i}\cr}}$}}&G^r&
\longrightarrow&0\cr
&&\Big\downarrow\rlap{$\vcenter{\hbox{$\scriptstyle
g$}}$}&&\Big\downarrow\rlap{$\vcenter{\hbox{$\scriptstyle
A$}}$}&&\Big\Vert&&\cr
0&\longrightarrow&G^s&
\longrightarrow&G^s\times G^r&
\longrightarrow&G^r&
\longrightarrow&0\cr}
$$
where $\pmatrix{\pi^{n_i}\cr}$ and $A$ denote the homomorphisms
$$\pmatrix{\pi^{n_1}&\cdots&0\cr \vdots&\ddots&\vdots\cr
0&\cdots&\pi^{n_r}\cr}\qquad\hbox{and}\qquad
\pmatrix{f_{11}&\cdots&f_{r1}\cr
\vdots&&\vdots\cr  f_{1s}&\cdots&f_{rs}\cr
\pi^{n_1}&\cdots&0\cr \vdots&\ddots&\vdots\cr0&\cdots&\pi^{n_r}\cr}
$$
respectively. The diagram 
has exact rows and columns
and, working in the category of fppf sheaves we deduce that
$$
K\cong K_1={\rm ker}(G^r\mathop{\longrightarrow}\limits^{A}G^{r+s}).
$$
Since the ring $R$ is a principal ideal domain, there exist an
invertible $r\times r$-matrix~$B$ and an invertible $(r+s)\times(r+s)$-matrix~$B'$
both with entries in $R$ so that 
$$
B'AB=\pmatrix{g_1&\cdots&0\cr\vdots&\ddots&\vdots\cr
0&\cdots&g_r\cr0&\cdots&0\cr\vdots&&\vdots\cr0&\cdots&0\cr}
$$
for certain $g_i\in R$. This shows that $K$ is isomorphic to the kernel of the map
$$
\pmatrix{g_1&\cdots&0\cr\vdots&\ddots&\vdots\cr
0&\cdots&g_r\cr}:G^r\longrightarrow G^r.
$$
This proves the corollary.

\proclaim Theorem 8.3. Let the ring $O$ and the category $\underline{D}$ be as above and let
$G=\{G_n\}$ be a $p$-divisible group over~$O$. Suppose that $R={\rm End}(G)$ is a discrete valuation
ring with uniformizer~$\pi$ and residue field $k=R/\pi R$. Suppose that
\smallskip
\item{--} every group scheme $G_n$ is an  object in the category~$\underline{D}$.
\item{--} the  map
$$
{\rm
Hom}_O(\gp,\gp)\mathop{\longrightarrow}\limits^{\delta} {\rm
Ext}^1_{\underline{D}}(\gp,\gp)
$$ 
associated to the exact sequence $0\rightarrow \gp\rightarrow G[\pi^2]\rightarrow\gp\rightarrow 0$ is
an isomorphism of 1-dimensional
$k$-vector spaces.
\smallskip
\noindent
Let $H=\{H_n\}$ be a $p$-divisible group over~$O$ for which the following hold.
\smallskip\item{--} every group scheme $H_n$ is an object in~$\underline{D}$;
\smallskip\item{--} each $H_n$ admits a filtration with flat closed subgroup schemes and successive
quotients isomorphic to~$G[\pi]$.
\smallskip\noindent
Then $H$ is isomorphic to $G^g$ for some $g\ge 0$.

\smallskip
\noindent{\bf Proof.} By Cor.~8.2. the group scheme $H[p^n]$ is for each $n\ge 1$ 
isomorphic to a group scheme of the shape $\mathop{\oplus}\limits_{i=1}^rG[\pi^{n_i}]$. 
Let $\overline{F}$ be an algebraic closure of the quotient field $F$ of~$O$. The
$\overline{F}$-points of $A[p^n]$ form a group isomorphic to $(\ZZ/p^n\ZZ)^g$ where $g={\rm dim}\,A$.
Therefore every direct summand of $A[p^n](\overline{F})$ is free over~$\ZZ/p^n\ZZ$. This implies 
that
$$
A[p^n]\cong\mathop\oplus\limits_{i=1}^rG[\pi^{en}]\cong G^r[p^n].
$$
Here $e$ denotes the ramification index of $R$ over~$\ZZ_p$ and~$re\cdot{\rm
dim}\,G[\pi](\overline{F})=g$.
Since the groups ${\rm Hom}(H[p^n], G^g[p^{n}])$ are finite, a compactness
argument shows that there is a cofinal projective system of such isomorphisms and hence
an isomorphism  $H\longrightarrow G^g$ as required.

\bibliography

\item{[\AB]} Abra\v skin, V.A.: Galois moduli of period $p$ group
schemes over a ring of Witt vectors, {\sl Izv. Ak. Nauk CCCP}, Ser.
Matem. {\bf 51}, (1987). English translation in {\sl Math. USSR
Izvestiya} {\bf 31} (1988) 1--46.
\item{[\ANT]}  Birch, B. and Kuyk, W. Eds.: {\sl Modular functions in one variable} IV.
Lecture Notes in Math. {\bf 476}, Springer-Verlag, New York 1975.
\item{[\BK]} Brumer, A. and Kramer, K.: Non-existence of certain semistable abelian
varieties, {\sl \hyphenation{Ma-nu-scrip-ta} Manuscripta Math.} {\bf 106} (2001) 291--304.
\item{[\FCALE]} Calegari, F.: Semistable abelian varieties over $\QQ$, preprint 2001. 
{\tt http://www.math.uiuc.edu /Algebraic-Number-Theory/0284/index.html}
\item{[\ED]} Edixhoven, S.J.: Minimal resolution and stable reduction of
$X_0(N)$,  Annales de l'Institut Fourier {\bf 40} (1990) 31--67. 
\item{[\FAL]} Faltings, G.: Endlichkeitss\"atze f\"ur abelsche
Variet\"aten \"uber Zahlk\"orpern, {\sl Invent. Math.} {\bf 73}
(1983) 349--366. 
\item{[\F]} Fontaine, J.-M.: Il n'y a pas de
vari\'et\'e ab\'elienne sur $\ZZ$, {\sl Invent. Math.} {\bf 81},
(1985) 515--538.
\item{[\SGA7]} Grothendieck, A.: Mod\`eles de N\'eron et monodromie, Exp IX in
{\sl Groupes de monodromie en g\'eom\'etrie alg\'ebrique}, SGA 7, Part I,  Lecture Notes in
Mathematics {\bf 288} (1971) Springer-Verlag, New York.
\item{[\KM]} Katz, N. and Mazur, B.: {\sl Arithmetic moduli of
elliptic curves}, Annals of Math. Studies {\bf 108}, Princeton
University Press, Princeton 1985.
\item{[\MA]} Martinet, J.: Petits discriminants des corps de nombres, 
in J.V.~Armitage, {\sl Journ\'ees Arithm\'et- iques 1980}, CUP
Lecture Notes Series {\bf 56}, Cambridge University Press, Cambridge 1981.
\item{[\MAZ]} Mazur, B.: Modular curves and the Eisenstein ideal, {\sl Publ. Math. IHES}, {\bf 47}
(1976), 33--186.
\item{[\ME]} Mestre, J.-F.:  Formules explicites et minorations de conducteurs de vari\'et\'es
alg\'ebriques, {\sl Compositio Math.} {\bf 58} (1986), 209--232.
\item{[\ODL]} Odlyzko, A.M.:
Unconditional bounds for discriminants, 1976.
{\tt http://www.dtc.umn.edu/}
{\tt $\sim$odlyzko/unpublished/discr.bound.table2}
\item{[\OPDV]} Oort, F.: Minimal $p$-divisible groups, Manuscript 3-IX-2002.
\item{[\RAY]} Raynaud, M.: Sch\'emas en groupes de type
$(p,\ldots,p)$, {\sl Bull. Soc. Math. France} {\bf 102}
(1974) 241--280.
\item{[\SCH]} Schoof, R.: Abelian varieties over cyclotomic fields with good reduction everywhere,
{\sl Math. Annalen.}, to appear.
\item{[\SCHO]} Schoof, R.: Abelian varieties over $\QQ(\sqrt{6})$ 
with good reduction everywhere, pp. 287--306 in: Miyake, K.: {\sl Class Field Theory ---Its Centenary
and Prospect},  Advanced Studies in Pure Mathematics {\bf 148}, Math. Soc. Japan, Tokyo 2001.
\item{[\TDB]} Tate J.T.: $p$-divisible groups, {\sl Proc. Conf. on local fields (Driebergen 1966)},
118--131. Springer-Verlag, Berlin.
\item{[\OT]} Tate, J.T. and Oort, F.: Group
schemes of prime order, {\sl Ann. Scient. \'Ecole Norm. Sup.} {\bf 3} (1970) 1--21.
\item{[\TAU]} Taussky, O.: A remark on the class field tower, {\sl J. London Math. Soc.} {\bf
12} (1937), 82--85.
\item{[\WAS]} Washington, L.C.: {\sl Introduction to cyclotomic fields},
Graduate Texts in Math. {\bf 83}, Springer-Verlag, Berlin Heidelberg
New York 1982.

\bye